\documentclass[12pt,british,reqno]{article}

\usepackage{babel}
\usepackage{amsmath,amssymb,amsthm}
\usepackage{ucs} 
\usepackage[utf8x]{inputenc} 
\usepackage[T1]{fontenc} 
\usepackage{eucal}
\usepackage{enumerate}
\usepackage[obeyspaces]{url}
\usepackage{calc}

\usepackage[matrix,arrow]{xy}
\usepackage[noBBpl]{mathpazo}
\usepackage{mathabx}

\makeatletter

\renewcommand{\p@enumii}{}

\def\@enum@{\list{\csname label\@enumctr\endcsname}%
           {\usecounter{\@enumctr}\def\makelabel##1{
\normalfont\ignorespaces\emph{{##1}~}}
\setlength{\labelsep}{3pt}
\setlength{\parsep}{0pt}
\setlength{\itemsep}{0pt}
\setlength{\leftmargin}{0pt}
\setlength{\labelwidth}{0pt}
\setlength{\listparindent}{\parindent}
\setlength{\itemsep}{0pt}
\setlength{\itemindent}{0pt}
\topsep=3pt plus 1pt minus 1 pt}}

\makeatother

\renewcommand{\epsilon}{\ensuremath{\varepsilon}}
\renewcommand{\phi}{\ensuremath{\varphi}}

\renewcommand{\to}{\ensuremath{\longrightarrow}}
\renewcommand{\mapsto}{\ensuremath{\longmapsto}}

\newcommand{\rp}{\ensuremath{\mathbb{R}P^2}}

\newcommand{\N}{\ensuremath{\mathbb N}}
\newcommand{\Z}{\ensuremath{\mathbb Z}}
\newcommand{\dt}{\ensuremath{\mathbb D}^{2}}
\newcommand{\St}[1][2]{\ensuremath{\mathbb S}^{#1}}
\newcommand{\FF}{\ensuremath{\mathbb F}}
\newcommand{\F}[1][n]{\ensuremath{\FF_{{#1}}}}
\newcommand{\sn}[1][n]{\ensuremath{S_{{#1}}}}
\newcommand{\an}[1][n]{\ensuremath{A_{{#1}}}}
\DeclareMathOperator{\id}{\text{Id}}

\renewcommand{\ker}[1]{\ensuremath{\operatorname{\text{Ker}}\left({#1}\right)}}
\newcommand{\im}[1]{\ensuremath{\operatorname{\text{Im}}\left({#1}\right)}}

\newcommand{\quat}[1][8]{\ensuremath{\mathcal{Q}_{#1}}}
\newcommand{\dih}[1]{\ensuremath{\operatorname{\text{Dih}}_{#1}}}
\newcommand{\dic}[1]{\ensuremath{\operatorname{\text{Dic}}_{#1}}}

\newcommand{\pnm}[1][n]{\ensuremath{P_{{#1}}(M)}}
\newcommand{\gpab}[1][G]{\ensuremath{{#1}\textsuperscript{Ab}}}
\newcommand{\bnab}[1]{\ensuremath{\left(B_{#1}(\rp)\right)\textsuperscript{Ab}}}

\makeatletter
\def\@map#1#2[#3]{\mbox{$#1 \colon\thinspace #2 \to #3$}}
\def\map#1#2{\@ifnextchar [{\@map{#1}{#2}}{\@map{#1}{#2}[#2]}}
\DeclareRobustCommand*\textsubscript[1]{\@textsubscript{\selectfont#1}}
\def\@textsubscript#1{{\m@th\ensuremath{_{\mbox{\fontsize\sf@size\z@#1}}}}}
\makeatother

\DeclareRobustCommand*{\up}[1]{\textsuperscript{#1}}

\newcommand{\ft}[1][n]{\ensuremath{\Delta_{#1}^2}}

\newcommand{\garside}[1][n]{\ensuremath{\Delta_{#1}}}

\newcommand{\brak}[1]{\ensuremath{\left\{ #1 \right\}}}
\newcommand{\ang}[1]{\ensuremath{\langle #1\rangle}}
\newcommand{\set}[2]{\ensuremath{\brak{#1 \,\mid\, #2}}}
\newcommand{\setang}[2]{\ensuremath{\ang{#1 \,\mid\, #2}}}

\newcommand{\setangl}[2]{\ensuremath{\ang{\left. #1 \,\right\rvert \, #2}}}

\newcommand{\setl}[2]{\ensuremath{\brak{\left. #1 \,\right\rvert \, #2}}}

\newcommand{\si}[2][{}]{\ensuremath{\sigma_{#2}^{#1}}}
\newcommand{\sii}[2][1]{\ensuremath{\sigma_{#2}^{-{#1}}}}

\theoremstyle{plain}
\newtheorem{thm}{Theorem}
\newtheorem{lem}[thm]{Lemma}
\newtheorem{prop}[thm]{Proposition}
\newtheorem{cor}[thm]{Corollary}

\theoremstyle{remark}
\newtheorem{rem}[thm]{Remark}
\newtheorem{rems}[thm]{Remarks}

\newcommand{\reth}[1]{Theorem~\protect\ref{th:#1}}
\newcommand{\relem}[1]{Lemma~\protect\ref{lem:#1}}
\newcommand{\repr}[1]{Proposition~\protect\ref{prop:#1}}
\newcommand{\reco}[1]{Corollary~\protect\ref{cor:#1}}
\newcommand{\resec}[1]{Section~\protect\ref{sec:#1}}
\newcommand{\rerem}[1]{Remark~\protect\ref{rem:#1}}

\newcommand{\req}[2][{}]{equation~(\protect\ref{eq:#2}\textsubscript{${#1}$})}
\newcommand{\reqref}[2][{}]{(\protect\ref{eq:#2}\textsubscript{${#1}$})}

\begin{document}

\title{The lower central and derived series of the braid groups
of the projective plane}

\author{Daciberg~Lima~Gon\c{c}alves\\
Departamento de Matem\'atica - IME-USP,\\
Caixa Postal~\textup{66281}~-~Ag.~Cidade de S\~ao Paulo,\\ CEP:~\textup{05314-970} - S\~ao Paulo - SP - Brazil.\\ 
e-mail: \url{dlgoncal@ime.usp.br}\vspace*{4mm}\\
John~Guaschi\\
Laboratoire de Math\'ematiques Nicolas Oresme UMR CNRS~\textup{6139},\\ Universit\'e de Caen BP~\textup{5186},
14032 Caen Cedex, France.\\
e-mail: \url{guaschi@math.unicaen.fr}}

\date{14\up{th} April 2010}

\begingroup
\renewcommand{\thefootnote}{}
\endgroup 

\maketitle


\begin{abstract}
\noindent
\emph{In this paper, we determine the lower central and
derived series for the braid groups of the projective plane. We are motivated in
part by the study of Fadell-Neuwirth short exact
sequences, but the problem is interesting in its own
right.}

\emph{The $n$-string braid groups $B_n(\rp)$ of the projective plane $\rp$ were originally studied by 
Van Buskirk during the 1960's, and are of particular interest due to the fact that they have torsion. The group $B_1(\rp)$ (resp.\ $B_2(\rp)$) is isomorphic to the cyclic group $\Z_{2}$ of order $2$ (resp.\ the generalised quaternion group of order $16$) and hence their lower central and derived series are known. If $n>2$, we first prove that the lower central series of $B_n(\rp)$ is constant from the commutator subgroup onwards. We observe that $\Gamma_2(B_3(\rp))$ is isomorphic to $(\F[3]\rtimes \quat)\rtimes \Z_{3}$, where $\F[k]$ denotes the free group of rank $k$, and $\quat$ denotes the quaternion group of order $8$, and that $\Gamma_2(B_4(\rp))$ is an extension of an index $2$ subgroup $K$ of $P_{4}(\rp)$ by $\Z_{2}\oplus \Z_{2}$. 
As for the derived series of
$B_n(\rp)$, we show that for all $n\geq 5$, it is constant from the
derived subgroup onwards. The group $B_n(\rp)$ being finite and
soluble for $n\leq 2$, the critical cases are $n=3,4$.
We are able to determine completely the derived
series of $B_3(\rp)$. The subgroups  $(B_3(\rp))^{(1)}$,  $(B_3(\rp))^{(2)}$ and  $(B_3(\rp))^{(3)}$ are isomorphic respectively to $(\F[3]\rtimes \quat)\rtimes \Z_{3}$, $\F[3]\rtimes \quat$ and $\F[9] \times \Z_2$, and we compute the derived series quotients of these groups. From $(B_3(\rp))^{(4)}$ onwards, the derived series of $B_{3}(\rp)$, as well as its successive derived series quotients, coincide
with those of $\F[9]$. We analyse the derived series of $B_{4}(\rp)$ and its quotients up to $(B_4(\rp))^{(4)}$, and we show that $(B_4(\rp))^{(4)}$ is a semi-direct of $\F[129]$ by $\F[17]$. Finally, we give a presentation of $\Gamma_2(B_n(\rp))$.}
\end{abstract}

\maketitle

\section{Introduction}

\subsection{Generalities and definitions}

Let $n\in\N$. The braid groups of the plane $\mathbb{E}^2$, denoted by $B_n$, and known as \emph{Artin braid groups}, were introduced by E.~Artin in~1925~\cite{A1,A2,A3}, and admit the following well-known presentation: $B_n$ is generated by elements $\sigma_1,\ldots, \sigma_{n-1}$, subject to the classical \emph{Artin relations}:
\begin{equation*}
\left\{
\begin{aligned}
&\text{$\si{i}\si{j}=\si{j}\si{i}$ if $\lvert i-j\rvert\geq 2$ and $1\leq i,j\leq n-1$}\\
&\text{$\si{i}\si{i+1}\si{i}=\si{i+1}\si{i}\si{i+1}$ for all $1\leq i\leq n-2$.}
\end{aligned}
\right.
\end{equation*}
A natural generalisation to braid groups of arbitrary topological
spaces was made at the beginning of the 1960's by Fox (using the
notion of configuration space)~\cite{FoN}. The braid groups of
compact, connected surfaces have been widely studied; (finite)
presentations were obtained in~\cite{Z1,Z2,Bi1,S}. As well as being
interesting in their own right, braid groups have played an important
r\^ole in many branches of mathematics, for example in topology,
geometry, algebra and dynamical systems, and notably in the study of
knots and links~\cite{BZ}, of the mapping class groups~\cite{Bi2,Bi3},
and of configuration spaces~\cite{CG,FH}. The reader may
consult~\cite{Bi2,Han,R} for some general references on the theory of
braid groups.

Let $M$ be a connected manifold of dimension~$2$ (or
\emph{surface}), perhaps with boundary. Further, we shall suppose that
$M$ is homeomorphic to a compact $2$-manifold with a finite (possibly
zero) number of points removed from its interior. We recall two
(equivalent) definitions of surface braid groups. The first is that
due to Fox. Let $F_n(M)$ denote the \emph{$n\up{th}$ configuration
space} of $M$, namely the set of all ordered $n$-tuples of distinct
points of $M$:
\begin{equation*}
F_n(M)=\set{(x_1,\ldots,x_n)}{\text{$x_i\in M$ and $x_i\neq x_j$ if $i\neq
j$}}. 
\end{equation*}
Since $F_n(M)$ is a subspace of the $n$-fold Cartesian product of $M$ with itself, the topology on
$M$ induces a topology on $F_n(M)$. Then we define the \emph{$n$-string pure braid group $\pnm$} of
$M$ to be $\pnm=\pi_1(F_n(M))$.
There is a natural action of the symmetric group $S_n$ on $F_n(M)$ by
permutation of coordinates, and the resulting orbit space $F_n(M)/\sn$
shall be denoted by $D_n(M)$. The fundamental group $\pi_1(D_n(M))$ is
called the \emph{$n$-string (full) braid group} of $M$, and shall be
denoted by $B_n(M)$. Notice that the projection $F_n(M) \to D_n(M)$ is
a regular $n!$-fold covering map. It is well known that $B_n$ is isomorphic to $B_n(\dt)$, and that the subgroup $P_n$ of pure braids of $B_{n}$ is isomorphic to $P_n(\dt)$, where $\dt$ is the closed $2$-disc.

The second definition of surface braid groups is geometric. Let
$\mathcal{P}=\brak{p_1, \dots, p_n}$ be a set of $n$ distinct points
of $M$. A \emph{geometric braid} of $M$ with basepoint $\mathcal{P}$
is a collection $\beta=(\beta_1, \dots, \beta_n)$ of $n$ paths
$\map{\beta}{[0,1]}[M]$ such that:
\begin{enumerate}[(a)]
\item for all $i=1,\ldots,n$, $\beta_i(0)=p_i$ and $\beta_i(1)\in \mathcal{P}$.
\item for all $i,j=1,\ldots,n$ and $i\neq j$, and for all $t\in [0,1]$, $\beta_i(t)\neq \beta_j(t)$.
\end{enumerate}
Two geometric braids are said to be \emph{equivalent} if there exists a homotopy between them through geometric braids. The usual concatenation of paths induces a group operation on the set of equivalence classes of geometric braids. This group is isomorphic to $B_n(M)$, and does not depend on the choice of $\mathcal{P}$. The subgroup of \emph{pure} braids, satisfying additionally $\beta_i(1)=p_i$ for all $i=1,\ldots,n$, is isomorphic to $P_n(M)$. There is a natural surjective homomorphism $\map{\tau}{B_n(M)}[\sn]$ which to a geometric braid $\beta$ associates the permutation $\tau(\beta)$ defined by $\beta_i(1)= p_{\tau(\beta)(i)}$. The kernel is precisely $P_n(M)$, and we thus obtain the following short exact sequence:
\begin{equation}\label{eq:perm}
1\to P_n(M)\to B_n(M) \stackrel{\tau}{\to} \sn\to 1.
\end{equation}

In this paper, we shall be primarily interested in the braid groups of
the real projective plane $\rp$. Along with the braid groups of the $2$-sphere, they are of particular interest, notably because they have non-trivial centre (which is also the case for the Artin braid groups), and torsion elements (which were characterised by
Murasugi~\cite{M}, see also~\cite{GG10}). We recall briefly some of their properties. If $\dt\subseteq \rp$ is a topological disc, there is a group homomorphism $\map{\iota}{B_n(\dt)}[B_n(\rp)]$ induced by the inclusion. If $\beta\in
B_n(\dt)$ then its image $\iota(\beta)$ shall be denoted simply by
$\beta$. A presentation of $B_{n}(\rp)$ was given in~\cite{vB} (see \repr{present}); in~\cite{GG4}, we obtained a presentation of $P_{n}(\rp)$. The first two
braid groups of $\rp$ are finite: $B_1(\rp)$ and $B_2(\rp)$ are isomorphic to $\Z_{2}$ and $\quat[16]$ respectively, where for $m\geq 2$, $\quat[4m]$ denotes the generalised quaternion group of order $4m$~\cite{vB}. If $n\geq 3$ then
$B_n(\rp)$ is infinite. For $n=3$, the Fadell-Neuwirth short exact sequence of pure braid groups yields the fact that $P_{3}(\rp)$ is isomorphic to a semi-direct product of a free group of rank two by $\quat$. If $n\geq 2$, the so-called `full twist' 
braid $\ft= (\sigma_1\cdots\sigma_{n-1})^n$ generates the centre
$Z(B_n(\rp))$ of $B_n(\rp)$, and is the unique element of $B_n(\rp)$
of order $2$. Here $\garside$ denotes the Garside (or `half twist') element of $B_{n}(\rp)$, defined by
\begin{equation*}
\garside = (\sigma_1 \cdots \sigma_{n-1}) (\sigma_1 \cdots \sigma_{n-2}) \cdots
(\sigma_1 \sigma_2)\sigma_1.
\end{equation*}
Further, the torsion of $B_{n}(\rp)$ is $4n$ and $4(n-1)$, and that of $P_{n}(\rp)$ is $2$ and $4$~\cite{GG3}. In~\cite{GG7}, we classified the virtually cyclic subgroups of $B_{n}(\rp)$ for all $n\in \N$, and in~\cite{GG10}, we characterised the finite subgroups of $B_{n}(\rp)$.

Our aim in this paper is to study the lower central and derived series
of the braid groups of $\rp$. We recall some definitions and
notation concerning these series. If $G$ is a group, then its
\emph{lower central series} $\brak{\Gamma_i(G)}_{i\in\N}$ is defined
inductively by $\Gamma_1(G)=G$, and $\Gamma_{i+1}(G)=[G, \Gamma_i(G)]$
for all $i\in\N$, and its \emph{derived series}
$\brak{G^{(i)}}_{i\in\N\cup\brak{0}}$ is defined inductively by
$G^{(0)}=G$, and $G^{(i)}=[G^{(i-1)}, G^{(i-1)}]$ for all $i\in\N$.
One may check easily that $\Gamma_i(G)\supseteq \Gamma_{i+1}(G)$ and
$G^{(i-1)}\supseteq G^{(i)}$ for all $i\in\N$, and for all $j\in\N$,
$j>i$, $\Gamma_j(G)$ (resp.\ $G^{(j)}$) is a normal subgroup of
$\Gamma_i(G)$ (resp.\ $G^{(i)}$).  Notice that $\Gamma_2(G)=G^{(1)}$
is the \emph{commutator subgroup} of $G$. The \emph{Abelianisation} of
the group $G$, denoted by $\gpab$ is the quotient $G/\Gamma_2(G)$; the
\emph{Abelianisation} of an element $g\in G$ is its
$\Gamma_2(G)$-coset in $\gpab$. The group $G$ is said to be
\emph{perfect} if $G=G^{(1)}$, or equivalently if $\gpab=\brak{1}$.
Following P.~Hall, for any group-theoretic property $\mathcal{P}$, a
group $G$ is said to be \emph{residually $\mathcal{P}$} if for any
(non-trivial) element $x \in G$, there exists a group $H$ with the
property  $\mathcal{P}$ and a surjective homomorphism
$\map{\phi}{G}[H]$ such that  $\phi(x) \neq 1$. It is well known that
a group $G$ is \emph{residually nilpotent} (respectively
\emph{residually soluble}) if and only if $\bigcap_{i \geq 1}\,
\Gamma_i(G)=\{ 1\}$ (respectively $\bigcap_{i \geq 0}\, G^{(i)}=\{
1\}$). If $g,h\in G$ then $[g,h]=gh g^{-1}h^{-1}$ will denote their
commutator.

The lower central series of groups and their successive quotients
$\Gamma_i/\Gamma_{i+1}$ are isomorphism invariants, and have been
widely studied using commutator calculus, in particular for free
groups of finite rank~\cite{H,MKS}. Falk and Randell, and
independently Kohno investigated the lower central series of the pure
braid group $P_n$, and were able to conclude that $P_n$ is residually
nilpotent~\cite{FR1,Ko}. Falk and Randell also studied the lower
central series of generalised pure braid groups~\cite{FR2,FR3}. 
Using the Reidemeister-Schreier rewriting process, Gorin and Lin
obtained a presentation of the commutator subgroup of $B_n$ for $n\geq
3$~\cite{GL}. For $n\geq 5$, they were able to
infer that $(B_n)^{(1)}=(B_n)^{(2)}$, and so $(B_n)^{(1)}$ is perfect.
From this it follows that $\Gamma_2(B_n)=\Gamma_3(B_n)$, hence $B_n$
is not residually nilpotent. If $n=3$ then they showed that
$(B_3)^{(1)}$ is a free group of rank~$2$, while if $n=4$, they proved
that $(B_4)^{(1)}$ is a semi-direct product of two free groups of
rank~$2$. By considering the action, one may see that
$(B_4)^{(1)}\supsetneqq (B_4)^{(2)}$. The work of Gorin and Lin on
these series was motivated by the study of almost periodic solutions
of algebraic equations with almost periodic coefficients. In~\cite{GG9,GG6} we studied the lower central and derived series of the $2$-sphere $\St$ and the finitely-punctured $2$-sphere. For $\St$, the case $n=4$ is critical, in the sense that if $n\neq 4$, $B_{n}(\St)$ is residually soluble if and only if $n<4$. It is an open question as to whether $B_{4}(\St)$ is residually soluble.

The above comments indicate that the study of the lower central and derived series of the braid groups of $\rp$ is an important problem in its own right, and it helps us to understand better the structure of such groups. But we are also motivated by the interesting question of the existence of a section (the \emph{`splitting problem')} for the following two short exact sequences of braid groups (notably for the case $M=\rp$) obtained by considering the long exact sequences in homotopy of fibrations of the corresponding configuration spaces:
\begin{enumerate}[(a)]
\item let $m,n\in\N$ and $m>n$. Then we have the \emph{Fadell-Neuwirth short exact sequence of pure
braid groups}~\cite{FaN}:
\begin{equation}\label{eq:split}
1 \to P_{n}(M\setminus\brak{x_1,\ldots,x_m}) \stackrel{i_{\ast}}{\to} P_{m+n}(\rp) \stackrel{p_{\ast}}{\to} P_{m}(\rp) \to 1,
\end{equation}
where $m\geq 3$ if $M=\St$~\cite{Fa,FvB}, $m\geq 2$ if $M=\rp$~\cite{vB}, and $m\geq 1$ otherwise~\cite{FaN}, and where $p_{\ast}$ is the group homomorphism which geometrically corresponds to forgetting the last $n$ strings, and $i_{\ast}$ is
inclusion (we consider $P_{n}(M\setminus\brak{x_1,\ldots,x_m})$ to be the subgroup of $P_{m+n}(\rp)$ of pure braids whose last $m$ strings are vertical). This short exact sequence plays a central r\^ole in the study of surface braid groups. It was used by~\cite{PR} to study mapping class groups, in the work of~\cite{GMP} on Vassiliev invariants for braid groups, as well as to obtain presentations for surface pure braid
groups~\cite{Bi1,S,GG1,GG4,GG6}. 
\item let $m,n\in\N$. Consider the group
homomorphism $\map{\tau}{B_{m+n}(M)}[{\sn[m+n]}]$, and let
$B_{m,n}(M)=\tau^{-1}(\sn[m]\times \sn[n])$ be the inverse image of
the subgroup $\sn[m]\times \sn[n]$ of $\sn[m+n]$. As in the pure braid
group case, we obtain a generalisation of the Fadell-Neuwirth
short exact sequence~\cite{GG2}:
\begin{equation}\label{eq:fnsurface}
1\to B_n(M\setminus\brak{x_1,\ldots,x_m}) \to B_{m,n}(M)\stackrel{p_{\ast}}{\to} B_m(M)\to 1,
\end{equation}
where we take $m\geq 3$ if $M=\St$, $m\geq 2$ if $M=\mathbb{R}P^2$ and $m\geq 1$ otherwise. Once more, $p_{\ast}$ corresponds geometrically to forgetting the last $n$ strings.
\end{enumerate}

We remark that if the above conditions on $n$ and $m$ are satisfied then the existence of a section for $p_{\ast}$ is equivalent to that of a geometric section for the corresponding configuration spaces (cf.~\cite{GG3,GG4}). The authors have recently solved the splitting problem for the short exact sequence~\reqref{split} for all surfaces~\cite{GG8}.
In~\cite{GG4}, we studied the short exact sequence~\reqref{fnsurface} in the case $M=\St$ of the sphere, and showed that if $m=3$ then~\reqref{fnsurface} splits if and only if $n\equiv 0,2\bmod 3$. Further, if $m\geq 4$ and~\reqref{fnsurface} splits then there exist $\epsilon_1,\epsilon_2\in \brak{0,1}$ such that $n\equiv \epsilon_1 (m-1)(m-2)-\epsilon_2 m(m-2) \bmod{m(m-1)(m-2)}$. An open question is whether this condition is also
sufficient.

Our main aim in this paper is to study the lower central and
derived series of the braid groups of $\rp$. This was motivated in part by the study of the problem of the existence of a section for the short exact sequences~\reqref{split}
and~\reqref{fnsurface}. To obtain a positive answer, it suffices of
course to exhibit an explicit section (although this may be easier
said than done!). However, and in spite of the fact that we possess
presentations of surface braid groups, in general it is very difficult
to prove directly that such an extension does not split. One of the
main methods that we used to prove the non-splitting of~\reqref{split}
for $n\geq 2$ and of~\reqref{fnsurface} for $m\geq 4$ was based on the
following observation: let $1\to K\to G\to Q\to 1$ be a split
extension of groups, where $K$ is a normal subgroup of $G$, and let
$H$ be a normal subgroup of $G$ contained in $K$. Then the extension
$1\to K/H\to G/H\to Q\to 1$ also splits. The condition on $H$ is
satisfied for example if $H$ is an element of either the lower central
series or the derived series of $K$. In~\cite{GG1}, considering the
extension~\reqref{split} with $n\geq 3$, we showed that it was
sufficient to take  $H=\Gamma_2(K)$ to prove the non-splitting of the
quotiented extension, and hence that of the full extension. In this
case, the kernel $K/\Gamma_2(K)$ is Abelian, which simplifies somewhat
the calculations in $G/H$. This was also the case in~\cite{GG4} for
the extension~\reqref{fnsurface} with $m\geq 4$. However, for the
extension~\reqref{split} with $n=2$, it was necessary to go a stage
further in the lower central series, and take $H=\Gamma_3(K)$. From
the point of view of the splitting problem, it is thus helpful to know
the lower central and derived series of the braid groups occurring in
these group extensions. But as we indicated earlier, these series are of course
interesting in their own right, and help us to understand better the structure of
surface braid groups. 

\subsection{Statement of the main results}

This paper is organised as follows. In \resec{lcsrp2}, we recall some general results concerning the splitting of the short exact sequence $1\to \Gamma_2(B_n(\rp)) \to B_n(\rp) \to  \bnab{n} \to 1$, where $\bnab{n}$ is the Abelianisation of $B_{n}(\rp)$, as
well as homological conditions for the stabilisation of the lower central series of a group (\relem{stallings}). We then go on to study the lower central series of $B_{n}(\rp)$, and we prove the following result.

\begin{thm}\label{th:lcsbn}
The lower central series of $B_n(\rp)$ is as follows.
\begin{enumerate}[(a)]
\item\label{it:lcsita} If $n=1$ then $B_1(\rp)=P_1(\rp)\cong \Z_{2}$, and $\Gamma_i(B_1(\rp))=\brak{1}$ for all $i\geq 2$.
\item\label{it:lcsitb} If $n=2$ then $B_2(\rp)$ is isomorphic to the generalised quaternion group $\quat[16]$ of order $16$. Its lower central series is given by $\Gamma_2(B_2(\rp))\cong \Z_4$, $\Gamma_3(B_2(\rp))\cong \Z_2$ and $\Gamma_i(B_2(\rp))=\brak{1}$  for all $i\geq 4$.
\item\label{it:lcsitc} For all $n\geq 3$, the lower central series of $B_n(\rp)$ is constant
from the commutator subgroup onwards: $\Gamma_m(B_n(\rp))=
\Gamma_2(B_n(\rp))$ for all $m\geq 2$. 
\end{enumerate}
Further, a presentation of $\Gamma_2(B_n(\rp))$ is given in \repr{fullpres}.
\end{thm}

The lower central series of $B_n(\rp)$ is thus completely determined.
In particular, for all $n\neq 2$, the lower central series of
$B_4(\rp)$ is constant from the commutator subgroup onwards, and $B_n(\rp)$ is residually nilpotent if and only if $n\leq 3$. A presentation of $\Gamma_{2}(B_{n}(\rp))$ is given in \repr{fullpres} in \resec{prescom}. The case $n=3$
is particularly interesting: as we shall see in \repr{gamma2rp34}, $\Gamma_2(B_3(\rp))$ is a semi-direct of the form $\left( \F[3]\rtimes \quat \right) \rtimes \Z_{3}$. This may be compared with Gorin and Lin's results for $\Gamma_2(B_3)$ and $\Gamma_2(B_4)$~\cite{GL} and with our result for $B_{4}(\St)$~\cite{GG9}. 

In \resec{dsrp2}, we study the derived series of $B_n(\rp)$. As in
the case of $B_n$ and $B_{n}(\St)$~\cite{GL,GG9}, $(B_n(\rp))^{(1)}$ is perfect if
$n\geq 5$, in other words, the derived series of $B_n(\rp)$ is
constant from $(B_n(\rp))^{(1)}$ onwards. The cases $n=1,2$ are
straightforward, and the groups $B_n(\rp)$ are finite and soluble. In
the case $n=3$, we make use of the semi-direct product decomposition
of $(B_3(\rp))^{(1)}$ of \repr{gamma2rp34}.

\begin{thm}\label{th:dsbn}
Let $n\in \N$, $n\neq 4$. The derived series of $B_n(\rp)$ is as follows.
\begin{enumerate}[(a)]
\item \label{it:ds1} If $n=1$ then $(B_n(\rp))^{(1)}=\brak{1}$.
\item \label{it:ds2} If $n=2$ then $(B_2(\rp))^{(1)}\cong\Z_4$ and $(B_2(\rp))^{(2)}=\brak{1}$.
\item  \label{it:ds3} Suppose that $n=3$. Then
\begin{enumerate}[(i)]
\item \label{it:ds3gam} $(B_3(\rp))^{(1)}=\Gamma_2(B_3(\rp))$ fits into the short exact sequence 
\begin{equation*}
1 \to  K \to  (B_3(\rp))^{(1)} \to \Z_3 \to 1,
\end{equation*}
where $K$ is an index $2$ subgroup of $P_3(\rp)$.
\item \label{it:ds3a} This short exact sequence splits; a section is given by associating $(\rho_{3}\sigma_{2} \sigma_{1})^4 \in (B_3(\rp))^{(1)}$ to a generator of $\Z_{3}$. The commutator subgroup $(B_3(\rp))^{(1)}$ is isomorphic to $(\F[3]\rtimes \quat)\rtimes \Z_{3}$, where the actions are given by \repr{gamma2rp34}.
\item \label{it:ds3b} We have $(B_3(\rp))^{(2)}\cong \F[3]\rtimes \quat$, where the action is given by \repr{gamma2rp34}. The quotient $(B_3(\rp))^{(1)}/(B_3(\rp))^{(2)}\cong \Z_3$, and there is a short exact sequence 
\begin{multline*}
1 \to  (B_3(\rp))^{(1)}/(B_3(\rp))^{(2)} \to  B_3(\rp)/(B_3(\rp))^{(2)} \to\\
B_3(\rp)/(B_3(\rp))^{(1)} \to 1,
\end{multline*}
where the extension $B_3(\rp))/(B_3(\rp))^{(2)}$ is isomorphic to the dihedral group $\dih{12}$ of order $12$.  Moreover, $(B_3(\rp))^{(2)}/(B_3(\rp))^{(3)}\cong\Z_{2}^4$, and $B_{3}(\rp)/(B_{3}(\rp))^{(3)}$ is an extension of $\Z_{2}^4$ by $\dih{12}$, so is of order $192$.
\item \label{it:ds3c} We have $(B_3(\rp))^{(3)} \cong \F[9]\oplus \Z_{2}$ and $(B_3(\rp))^{(3)}/(B_3(\rp))^{(4)} \cong \Z^9\oplus \Z_{2}$. Further, $B_3(\rp)/(B_3(\rp))^{(4)}$ is an extension of $\Z^9\oplus \Z_{2}$ by $B_{3}(\rp)/(B_{3}(\rp))^{(3)}$, so is infinite, and for all $i\geq 4$, $(B_3(\rp))^{(i)}\cong (\F[9])^{(i-3)}$.
\end{enumerate}
\item \label{it:ds4} If $n\geq 5$ then $(B_n(\rp))^{(2)}=(B_n(\rp))^{(1)}$, so $(B_n(\rp))^{(1)}$ is perfect. A presentation of $(B_n(\rp))^{(1)}$ is given in \repr{fullpres}.
\end{enumerate}
\end{thm}

So if $n\neq 4$, the derived series of $B_n(\rp)$ is thus completely determined (up to knowing the derived series of the free group $\F[9]$ of rank~$9$). In particular, if $n\neq 4$, $B_n(\rp)$ is residually soluble if and only if $n< 4$ (\reco{residb4}). We remark that part~(\ref{it:lcsitc}) of \reth{lcsbn} and the first statement of part~(\ref{it:ds4}) of \reth{dsbn} appeared in~\cite{BM} where the authors asserted that the results may be proved along the lines of our proof in the case of the sphere~\cite{GG9}. We give the details of the proofs. As for $B_{n}$ and $B_{n}(\St)$~\cite{GL,GG9,GG6}, the case $n=4$ is somewhat delicate. We are able to determine some of the terms and quotients of the derived series of $B_4(\rp)$.

\begin{thm}\label{th:dsb4}
Suppose that $n=4$. 
\begin{enumerate}[(a)]
\item \label{it:ds4gam} The group $(B_4(\rp))^{(1)}=\Gamma_2(B_4(\rp))$ is given by an extension 
\begin{equation*}
1 \to  K \to  (B_4(\rp))^{(1)} \to A_{4}  \to 1
\end{equation*}
where $K$ is a subgroup of $P_4(\rp)$ of index two.

\item \label{it:ds4a} 
\begin{enumerate}
\item\label{it:ds4ai} We have the following isomorphism:
\begin{equation*}
(B_4(\rp))^{(1)} \cong (B_4(\rp))^{(2)} \rtimes \Z_{3},
\end{equation*}
where the action on $(B_4(\rp))^{(2)}$ is given by conjugation by $(\rho_{3}\sigma_{2} \sigma_{1})^4$,
and
\begin{equation*}
(B_4(\rp))^{(1)}/(B_4(\rp))^{(2)} \cong \Z_{3}.
\end{equation*}

\item\label{it:ds4aii} We have a short exact sequence 
\begin{multline*}
1 \to  (B_4(\rp))^{(1)}/(B_4(\rp))^{(2)} \to  B_4(\rp)/(B_4(\rp))^{(2)} \to\\  
B_4(\rp)/(B_4(\rp))^{(1)} \to 1,
\end{multline*}
where $B_4(\rp)/(B_4(\rp))^{(2)}$ is isomorphic to the dihedral group $\dih{12}$ of order $12$. 

\item\label{it:ds4aiii} The group $(B_4(\rp))^{(2)}$ is given by an extension 
\begin{equation*}
1\to K\to (B_4(\rp))^{(2)} \to \Z_2\oplus\Z_2 \to 1.
\end{equation*}
\end{enumerate}

\item \label{it:ds4b} $(B_4(\rp))^{(2)}/(B_4(\rp))^{(3)}\cong \Z_{2}^4$, and $(B_4(\rp))^{(1)}/(B_4(\rp))^{(3)}\cong \Z_{2}^4\rtimes \Z_{3}$, where the action of $\Z_{3}$ permutes cyclically the three non-trivial elements of the first and second (resp.\ the third and fourth) copies of $\Z_{2}$. 
\item \label{it:ds4d} The group $(B_4(\rp))^{(3)}$ is a subgroup of $K$ of index four. Further,
\begin{equation*}
(B_4(\rp))^{(3)}\cong  (\F[5]\rtimes \F[3])\rtimes \Z_{4},
\end{equation*}
where the action is described by equations~\reqref{actz4}--\reqref{actf3f5c}. Moreover,
\begin{equation*}
(B_4(\rp))^{(3)}/(B_4(\rp))^{(4)} \cong \Z_{2}^8 \oplus \Z_{4},
\end{equation*}
and $(B_4(\rp))^{(4)}$ is a semi-direct product of the form $\F[129]\rtimes \F[17]$ where the action is that induced by $\F[3]$ on $\F[5]$. From $i=4$ onwards, we have $(B_4(\rp))^{(i+4)}\cong(\F[129]\rtimes \F[17])^{(i)}$ for all $i\geq 0$. 
\end{enumerate}
\end{thm}

A presentation of $B_4(\rp))^{(1)}$ derived from that of \repr{fullpres} is given during the proof of \reth{dsb4}. As in the case of $B_{4}(\St)$, it is an open question as to whether $B_4(\rp)$ is residually soluble or not. 

In~\cite{BGG}, the lower central series of braid groups of orientable
surfaces of genus $g\geq 1$, with and without boundary, was analysed. The study of the lower central series of non-orientable surfaces of genus at least two is the subject of work in progress.



\subsubsection*{Acknowledgements}

This work took place during the visit of the second author to the
Departmento de Matem\'atica do IME-Universidade de S\~ao Paulo during
the periods~14\up{th}--29\up{th}~April~2008, 18\up{th}~July--8\up{th} August~2008, 31\up{st}~October--10\up{th}~November~2008 and 17\up{st}~September--5\up{th}~October~2009, and of the visit of the first author to the Laboratoire de Math\'ematiques Nicolas Oresme, Universit\'e de Caen during the period 21\up{st}~November--21\up{st}~December~2008 and 19\up{th}~November--16\up{th}~December~2009. This work was supported by the international Cooperation USP/Cofecub project n\up{o} 105/06, by the CNRS/CNPq project n\up{o}~21119 and by the~ANR project TheoGar n\up{o} ANR-08-BLAN-0269-02.

\section{The lower central series of $B_n(\rp)$}\label{sec:lcsrp2}

The main aim of this section is to prove \reth{lcsbn}, which describes the lower central series of $B_n(\rp)$. Before doing so, we state some general results concerning $B_n(\rp)$, as well as some general homological conditions for the stabilisation of the lower central series of a group (\relem{stallings}). We start by recalling Van Buskirk's presentation of $B_{n}(\rp)$.


\begin{prop}[Van Buskirk~\cite{vB}]\label{prop:present}
Let $n\in\N$. The following constitutes a presentation of the group $B_n(\rp)$:\\
\underline{\textbf{generators:}}
$\si{1},\ldots,\si{n-1},\rho_{1},\ldots,\rho_{n}$.\\ 
\underline{\textbf{relations:}}
\begin{align*}
\si{i}\si{j} &=\si{j}\si{i}\quad\text{if $\lvert i-j\rvert\geq 2$,}\\
\si{i}\si{i+1}\si{i}&=\si{i+1}\si{i}\si{i+1} \quad\text{for $1\leq i\leq n-2$,}\\
\si{i}\rho_{j}&=\rho_{j}\si{i}\quad\text{for $j\neq i,i+1$,}\\
\rho_{i+1}&=\sii{i}\rho_{i}\sii{i} \quad\text{for $1\leq i\leq n-1$,}\\
\rho_{i+1}^{-1}\rho_{i}^{-1}\rho_{i+1}\rho_{i}&= {\si{i}}^2 \quad\text{for $1\leq i\leq n-1$,}\\
\rho_{1}^2&=\si{1}\si{2}\cdots\si{n-2}\si[2]{n-1} \si{n-2}\cdots\si{2}\si{1}.
\end{align*}
\end{prop}

\begin{rem}\label{rem:genspn}
Let $n\in \N$. It is well known that $\setl{B_{i,j},\, \rho_{k}}{1\leq i<j\leq n,\, 1\leq k\leq n}$ is a generating set for $P_{n}(\rp)$, where 
\begin{equation*}
B_{i,j}=\sigma_{j-1}\cdots \sigma_{i+1}\sigma_{i}^2 \sigma_{i+1}^{-1} \cdots \sigma_{j-1}^{-1}.
\end{equation*}
\end{rem}

Let $n\in\N$, let $\bnab{n}= B_n(\rp)/\Gamma_2(B_n(\rp))$ denote the
Abelianisation of $B_n(\rp)$, and let
$\map{\alpha}{B_n(\rp)}[\bnab{n}]$ be the canonical projection.  Then
we have the following short exact sequence:
\begin{equation}\label{eq:ses}
\xymatrix{%
1\ar[r] & \Gamma_2(B_n(\rp)) \ar[r] & B_n(\rp) \ar[r]^-{\alpha} & \bnab{n} \ar[r] & 1.}
\end{equation} 

We first prove the following result which deals with  this short exact sequence. 
\begin{prop}\label{prop:bnabz}
Let $n\in\N$. Then $\bnab{n}=B_n(\rp)/\Gamma_2(B_n(\rp)) \cong \Z_{2}\oplus \Z_{2}$, where the generators of the first (resp.\ second) copy of $\Z_{2}$ is the image 
of the generators $\sigma_i$ (resp.\ $\rho_j$).
\end{prop}

\begin{proof}
This follows easily by Abelianising the presentation of $B_n(\rp)$ given in \repr{present}. The generators $\sigma_i$ (resp.\ $\rho_j$) of $B_n(\rp)$ are all identified by $\alpha$ to a single generator $\overline{\sigma}=\alpha(\sigma_i)$ (resp.\ $\overline{\rho}=\alpha(\rho_j)$) of the first (resp.\ second) $\Z_{2}$-summand.
\end{proof}


We recall the following lemma from \cite{GG9}.

\begin{lem}[\cite{GG9}]\label{lem:stallings}
Let $G$ be a group, and let $\map{\delta}{H_2(G,\Z)}[H_2(\gpab,\Z)]$.
denote the homomorphism induced by Abelianisation. Then $\Gamma_2(G)=\Gamma_3(G)$ if and only if $\delta$ is surjective.
\end{lem}

We now come to the proof of \reth{lcsbn}.

\begin{proof}[Proof of \reth{lcsbn}]
Since $B_1(\rp)= \pi_{1}(\rp)$ and $B_2(\rp)\cong \quat[16]$\cite{vB}, parts~(\ref{it:lcsita}) and~(\ref{it:lcsitb}) follow easily. Now suppose that $n\geq 3$. First observe that $H_2(\Z_2\oplus \Z_2)\cong \Z_2$.  By \relem{stallings}, if the homomorphism $\delta$ is surjective then $\Gamma_2(B_n(\rp))=\Gamma_3(B_n(\rp))$, and part~(\ref{it:lcsitc}) follows. Otherwise, if $\delta$ is not surjective then it is trivial, and we obtain the following exact sequence:
\begin{equation*}
1\to \Z_2 \to \Gamma_2(B_n(\rp))/\Gamma_3(B_n(\rp)) \to H_1(B_n(\rp),\Z) \to \gpab[(B_n(\rp))]
\to 1.
\end{equation*} 
It follows that $\Z_2 \to \Gamma_2(B_n(\rp))/\Gamma_3(B_n(\rp))$ is an isomorphism. So we have the  short exact sequence:
\begin{equation*}
1 \to \Z_2 \to B_n(\rp)/\Gamma_3(B_n(\rp)) \to  \underbrace{H_1(B_n(\rp),\Z)}_{\Z_2 \oplus \Z_2}  \to  1,
\end{equation*} 
and hence the middle group, which we denote by $H$, is of order $8$. Since the quotient $\Gamma_2(B_n(\rp))/\Gamma_3(B_n(\rp))$ is non trivial, we conclude that $H$ is non Abelian, and so is either $\quat$ or the dihedral group $\dih{8}$.

We claim that there is no surjective homomorphism $B_n(\rp)\to H$. To see this, let $\map{\varphi}{B_n(\rp)}[H]$ be a homomorphism. Since $\si{i}\si{i+1}\si{i}=\si{i+1}\si{i}\si{i+1}$ for all $1\leq i\leq n-2$, the $\sigma_{i}$ are pairwise conjugate. Hence $\varphi(\sigma_i)$ and $\varphi(\sigma_j)$ are conjugate in $H$ for all $1\leq i,j\leq n-1$. But in both $\quat$ and $\dih{8}$, any two conjugate elements commute. Applying $\varphi$ to the relation $\si{i}\si{i+1}\si{i}=\si{i+1}\si{i}\si{i+1}$ and using induction yields $\varphi(\sigma_i)= \varphi(\sigma_j)$ for all $1\leq i,j\leq n-1$. If $\varphi(\sigma_i)=1$ then the relation $\rho_{i+1}=\sigma_{i}^{-1} \rho_{i} \sigma_{i}^{-1}$ implies that $\varphi(\rho_i)=\varphi(\rho_{i+1})$ for all $1\leq i\leq n-1$, and thus $\im{\varphi}=\ang{\varphi(\rho_{1})}\neq H$. So we may assume that $\varphi(\sigma_i)\neq 1$.

Suppose first that $n\geq 4$. Given $1\leq i\leq n$, there exists $1\leq j\leq n-1$ such that $\sigma_j$ commutes with $\rho_i$, and so $\varphi(\sigma_{j})$ commutes with $\varphi(\rho_i)$ for all $i$. If $n=3$ then a similar analysis shows that $\varphi(\rho_{1})$ and $\varphi(\rho_{3})$ commute with the $\varphi(\sigma_j)$. Further, $\rho_{2}=\sigma_{1}^{-1} \rho_{1} \sigma_{1}^{-1}$, and hence $\varphi(\rho_{2})$ commutes with the $\varphi(\sigma_j)$. In both cases, we conclude that $\im{\varphi}$ is contained in the centraliser of $\varphi(\sigma_{1})$ in $H$. A necessary condition for $\varphi$ to be surjective is that $\varphi(\sigma_{1})$ be central in $H$, and so $\varphi(\sigma_{1})$ must be of order $2$. Once more the relation $\rho_{i+1}=\sigma_{i}^{-1} \rho_{i} \sigma_{i}^{-1}$ implies that $\varphi(\rho_i)=\varphi(\rho_{i+1})$ for all $1\leq i\leq n-1$, and hence $\im{\varphi}=\ang{\varphi(\rho_1),\varphi(\sigma_{1})}\neq H$.

Thus no homomorphism $B_n(\rp)\to H$ is surjective, but this contradicts the surjectivity of the canonical projection $B_n(\rp) \to B_n(\rp)/\Gamma_3(B_n(\rp))$. This completes the proof of part~(\ref{it:lcsitc}), and thus that of \reth{lcsbn}.
\end{proof}

We may obtain a much better description of $\Gamma_{2}(B_{3}(\rp))$ as follows. This will be helpful in the analysis of the derived series in \resec{dsrp2}.

\begin{prop}\label{prop:gamma2rp34}
The group $\Gamma_2(B_3(\rp))$ is isomorphic to $\left( \F[3]\rtimes \quat \right) \rtimes \Z_{3}$.
The actions may be described as follows. Writing $\quat=\setang{x,y}{x^2=y^2,\, yxy^{-1}=x^{-1}}$, $\F[3]=\F[3](z_{1}, z_{2}, z_{3})$ and $\Z_{3}=\ang{u}$, we have:
\begin{align*}
xz_{1}x^{-1}&=z_{1}^{-1}  & xz_{2}x^{-1}&=z_{1}^{-1}z_{3}^{-1}z_{1} & xz_{3}x^{-1}&=z_{1}^{-1}z_{2}^{-1}z_{1}\\
yz_{1}y^{-1}&= z_{2}z_{3}z_{1} & yz_{2}y^{-1}&= z_{2}^{-1} & yz_{3}y^{-1}&= z_{2}z_{3}^{-1}z_{2}^{-1}\\
u z_{1}u^{-1} &=x^2  z_{3}z_{1} & u z_{2}u^{-1}&=x^2 z_{1}^{-1} & u z_{3}u^{-1}&=x^2 z_{2}^{-1}z_{1}^{-1} z_{3}^{-1}\\
u x u^{-1}&=xy & u y u^{-1}&=x, &
\end{align*}
where $u=(\rho_{3} \sigma_{2}\sigma_{1})^4$, $x=\rho_{2}\rho_{1}$, $y=\rho_{2} B_{1,2} \rho_{3}^{-1}$, $z_{1}=\rho_{3}^2$, $z_{2}=B_{2,3}$ and $z_{3}=\rho_{3}B_{2,3}\rho_{3}^{-1}$.
\end{prop}

\begin{rems}\label{rem:defab}\mbox{}
\begin{enumerate}[(a)]
\item The commutator subgroup of $B_{4}(\rp)$ will be analysed in more detail in \resec{dsrp2}.
\item\label{it:defab} Let $n\geq 2$. Recall from~\cite[Proposition~26]{GG3} that there exist two elements $a,b\in B_{n}(\rp)$ defined by:
\begin{equation}\label{eq:defab}
\left\{ \begin{aligned}
a &=\rho_{n} \sigma_{n-1}\cdots\sigma_{1}=\sigma_{n-1}^{-1} \cdots \sigma_{1}^{-1}\rho_{1}\\
b &=\rho_{n-1} \sigma_{n-2}\cdots\sigma_{1}=\sigma_{n-2}^{-1} \cdots \sigma_{1}^{-1}\rho_{1},
\end{aligned}\right.
\end{equation} 
of order $4n$ and $4(n-1)$ respectively. These elements satisfy~\cite[Remark~27]{GG3}:
\begin{equation*}
\text{$b^{n-1}=\rho_{n-1}\cdots \rho_{1}$ and $a^n=\rho_{n}\cdots\rho_{1}$.}
\end{equation*}
From~\cite[page~777]{GG3}, conjugation by $a^{-1}$ permutes cyclically the following two collections of elements:
\begin{equation}\label{eq:cyclicperm}
\left\{
\begin{gathered}
\sigma_{1},\ldots,\sigma_{n-1},\, a^{-1}\sigma_{n-1}a,\, \sigma_{1}^{-1},\ldots,\sigma_{n-1}^{-1},\, a^{-1}\sigma_{n-1}^{-1}a, \quad\text{and}\\
\rho_{1},\ \ldots,\rho_{n},\, \rho_{1}^{-1},\ldots, \rho_{n}^{-1}.
\end{gathered}
\right.
\end{equation}
In particular,
\begin{equation}\label{eq:cyclicperm2}
\left\{
\begin{aligned}
a^n \sigma_{i} a^{-n}&=\sigma_{i}^{-1} \quad \text{for all $1\leq i \leq n-1$}\\
a^n \rho_{j} a^{-n}&=\rho_{j}^{-1} \quad \text{for all $1\leq j \leq n$.}
\end{aligned}
\right.
\end{equation}
Further, for all $1\leq i\leq n$~\cite{GG10},
\begin{equation*}
\garside \rho_{i} \garside^{-1} =\rho_{n+1-i}^{-1}
\end{equation*}
in $B_{n}(\rp)$, which implies that
\begin{equation}\label{eq:conjgarside}
\garside a \garside^{-1}=\garside \rho_{n} \sigma_{n-1} \cdots \sigma_{1} \garside^{-1}=\rho_{1}^{-1} \sigma_{1} \cdots\sigma_{n-1}=a^{-1}.
\end{equation}
These observations will be used frequently in what follows.
\end{enumerate}
\end{rems}

\begin{proof}[Proof of \repr{gamma2rp34}]
Let $n\geq 2$, and let $\alpha$ be the Abelianisation homomorphism of \req{ses}, where $\alpha(\sigma_{i})=\overline{\sigma}$ and $\alpha(\rho_{j})=\overline{\rho}$. The permutation homomorphism $\tau$ of \req{perm} induces a homomorphism $\map{\overline{\tau}}{\bnab{n}}[\ang{\overline{\sigma}}]$, and we obtain the following commutative diagram of short exact sequences:
\begin{equation}\label{eq:diaggam2}
\begin{xy}*!C\xybox{%
\xymatrix{%
& 1 \ar[d] & 1 \ar[d] & 1 \ar[d] &\\
1 \ar[r] & K \ar[d]\ar[r] & \Gamma_{2}(B_{n}(\rp)) \ar[d]\ar[r]^{\tau'} & \an=\Gamma_{2}(\sn) \ar[d]\ar[r] & 1\\
1 \ar[r] & P_{n}(\rp) \ar[d]^{\alpha'}\ar[r] & B_{n}(\rp) \ar[d]^{\alpha}\ar[r]^{\tau} & \sn \ar[d]^h\ar[r] & 1\\
1 \ar[r] & \ang{\overline{\rho}} \ar[d]\ar[r] & \bnab{n} \ar[d]\ar[r]^{\overline{\tau}} & \ang{\overline{\sigma}} \ar[d]\ar[r] & 1.\\
& 1 & 1 & 1 & }}
\end{xy}
\end{equation}
Here $\tau'$ (resp. $\alpha'$) is the restriction of $\tau$ (resp.\ $\alpha$) to $\Gamma_{2}(B_{n}(\rp))$ (resp.\ to $P_{n}(\rp)$), $h$ is the homomorphism that to a transposition associates $\overline{\sigma}$, and $K=\ker{\alpha'}=\ker{\tau'}$ is of index $2$ in $P_{n}(\rp)$ (recall from \repr{bnabz} that $\bnab{n}\cong \Z_{2}\oplus \Z_{2}$, and $\ang{\overline{\rho}}\cong \ang{\overline{\sigma}}\cong \Z_{2}$).

Now let $n=3$. From~\cite{vB}, we know that
\begin{equation}\label{eq:vanbuskirk}
P_{3}(\rp)\cong \F[2]\rtimes \quat.
\end{equation}
Let us first determine generating sets of the two factors in terms of Van Buskirk's generators (this action was previously described in~\cite{GG4}, but in terms of a different generating set). From the Fadell-Neuwirth short exact sequence~\reqref{split}, we have
\begin{equation*}
1 \to \pi_{1}(\rp\setminus \brak{x_{1},x_{2}}) \to P_{3}(\rp) \stackrel{p_{\ast}}{\to} P_{2}(\rp)\cong \quat \to 1,
\end{equation*}
where $\pi_{1}(\rp\setminus \brak{x_{1},x_{2}})\cong \F[2]$ is a free group of rank two with basis $(\rho_{3},B_{2,3})$. The two elements $a=\rho_{3} \sigma_{2}\sigma_{1}$ and $b=\rho_{2}\sigma_{1}$ of \req{defab} are of order $12$ and $8$ respectively, and satisfy:
\begin{equation}\label{eq:a3b2}
\text{$b^2=\rho_{2}\rho_{1}$ and $a^3=\rho_{3}\rho_{2}\rho_{1}$.}
\end{equation}
From~\cite[Proposition~15]{GG10}, there is a copy of $\quat[16]$ in $B_{3}(\rp)$ of the form $\ang{b,\garside[3]a^{-1}}$, and by general arguments, one sees that it has two subgroups isomorphic to $\quat$, of the form $\ang{b^2,\garside[3]a^{-1}}$ and $\ang{b^2,b\garside[3]a^{-1}}$ respectively. We shall be interested in the latter copy since it is a subgroup of $P_{3}(\rp)$.

We have that $a^4$ is of order $3$, $a^4\in \ker{\alpha}=\Gamma_{2}(B_{3}(\rp))$, and $\tau(a^4)=\tau(a)=(1,2,3)$. Since $\an[3]=\ang{(1,2,3)}$, the correspondence $(1,2,3)\mapsto a^4$ defines a section for $\tau'$, and hence
\begin{equation}\label{eq:semigamb3}
\Gamma_{2}(B_{3}(\rp))\cong K\rtimes \Z_{3}
\end{equation}
from \req{diaggam2}. Let us now study the structure of $K$ in order to calculate the action.

By construction, $K$ is the kernel of $\alpha'$, and so is an index $2$ subgroup of $P_{3}(\rp)\cong \F[2]\rtimes \quat$. The homomorphism $\alpha'$ is defined on the generators of $P_{3}(\rp)$ (cf.\ \rerem{genspn}) by $\rho_{j}\mapsto \overline{\rho}$ for $j=1,2,3$, and for $1\leq i<j\leq 3$, $B_{i,j}$ is sent to the trivial element of $\ang{\overline{\rho}}$. Since
\begin{equation}\label{eq:bgara}
b\garside[3]a^{-1}=\sigma_{1}^{-1}\rho_{1} \ldotp \sigma_{1} \sigma_{2} \sigma_{1} \ldotp \sigma_{1}^{-1}\sigma_{2}^{-1}\rho_{3}^{-1}=\sigma_{1}^{-1}\rho_{1}\sigma_{1}\rho_{3}^{-1}=\rho_{2} \sigma_{1}^2 \rho_{3}^{-1}=\rho_{2} B_{1,2} \rho_{3}^{-1},
\end{equation}
we see that our copy $\ang{b^2,b\garside[3]a^{-1}}$ of $\quat$ lies in $\ker{\alpha'}$. Thus we have a commutative diagram of short exact sequences:
\begin{equation*}
\xymatrix{%
& 1\ar[d] & 1\ar[d] &  & \\
1 \ar[r]  & \ker{p_{\ast}\left\lvert_{K}\right.} \ar[r] \ar[d] & K \ar[r]^{p_{\ast}\left\lvert_{K}\right.} \ar[d] & P_{2}(\rp)\ar[r] \ar@{=}[d] & 1\\
1 \ar[r]  & \F[2](\rho_{3},B_{2,3}) \ar[r] \ar[d]^{\alpha'\left\lvert_{\F[2](\rho_{3},B_{2,3})}\right.} & P_{3}(\rp) \ar[r]^{p_{\ast}} \ar[d]^{\alpha'} & P_{2}(\rp) \ar[r]  & 1\\
1 \ar[r]  & \ang{\overline{\rho}} \ar@{=}[r] \ar[d] & \ang{\overline{\rho}}.  \ar[d] &  & \\
& 1 & 1 & & }
\end{equation*}
Note that we have used the following facts in order to construct this diagram:
\begin{enumerate}[(i)]
\item $\alpha'\left\lvert_{\pi_{1}(\rp\setminus \brak{x_{1},x_{2}}}\right.$ is surjective onto $\ang{\overline{\rho}}$ since $\rho_{3}\in \F[2](\rho_{3},B_{2,3})$.
\item $p_{\ast}(b^2)=(\rho_{2}\rho_{1})^2=\ft[2]$ is equal to $B_{1,2}$ in $P_{2}(\rp)$, and $p_{\ast}(b\garside[3]a^{-1})=\rho_{2} B_{1,2}$ by \req{bgara}, and we conclude that $p_{\ast}\left\lvert_{K}\right.$ is surjective onto $P_{2}(\rp)$.
\end{enumerate}
This second fact also implies that
\begin{equation}\label{eq:kkprime}
K\cong K'\rtimes \quat,
\end{equation}
where $K'=\ker{p_{\ast}\left\lvert_{K}\right.}=\ker{\alpha'}\cap \F[2](\rho_{3},B_{2,3})=\ker{\alpha'\left\vert_{\F[2](\rho_{3},B_{2,3})}\right.}$ is of index two in $\F[2](\rho_{3},B_{2,3})$. The application of the Reidemeister-Schreier rewriting process with Schreier transversal $\brak{1,\rho_{3}}$ to this restriction shows that $K'=\F[3]\left(\rho_{3}^2,B_{2,3}, \rho_{3}B_{2,3}\rho_{3}^{-1}\right)$ is a free group of rank $3$. Combining equations~\reqref{semigamb3} and~\reqref{kkprime}, we obtain $K\cong (\F[3]\rtimes \quat)\rtimes \Z_{3}$. The actions may be deduced from the action of $\quat=\ang{b^2,b\garside[3]a^{-1}}$ on $\F[2](\rho_{3},B_{2,3})$ which we now determine. \rerem{defab}(\ref{it:defab}) and~\req{a3b2} imply that conjugation by $b^2=\rho_{2}\rho_{1}$ is given by:
\begin{equation}\label{eq:actr12}
\left\{
\begin{aligned}
\rho_{2}\rho_{1} \ldotp \rho_{3} \ldotp \rho_{1}^{-1}\rho_{2}^{-1} &= \rho_{3}^{-1} a^3 \rho_{3} a^{-3} \rho_{3}=\rho_{3}^{-1}\\
\rho_{2}\rho_{1} \ldotp B_{2,3} \ldotp \rho_{1}^{-1}\rho_{2}^{-1} &= \rho_{3}^{-1} a^3 B_{2,3} a^{-3} \rho_{3}=\rho_{3}^{-1} B_{2,3}^{-1} \rho_{3}.
\end{aligned}\right.
\end{equation}
From this, we deduce that under conjugation by $b^2$,
\begin{equation}\label{eq:actconj1}
\left\{\begin{aligned}
\rho_{3} & \mapsto \rho_{3}^{-1}\\
B_{2,3}\rho_{3} & \mapsto (B_{2,3} \rho_{3})^{-1}.
\end{aligned}\right.
\end{equation}
As for conjugation by $b\garside[3]a^{-1}=\rho_{2} B_{1,2} \rho_{3}^{-1}$, we have
\begin{align}
\rho_{2} B_{1,2} \rho_{3}^{-1}\ldotp \rho_{3} \ldotp \rho_{3} B_{1,2}^{-1} \rho_{2}^{-1}&=\rho_{2} \rho_{3} \rho_{2}^{-1}= \rho_{3} \rho_{2} B_{2,3}^{-1}  \rho_{2}^{-1} \quad\text{using \repr{present}}\notag\\
& =\rho_{3} \rho_{2} \rho_{1} B_{2,3}^{-1} \rho_{1}^{-1} \rho_{2}^{-1} \rho_{3}^{-1} \rho_{3}= a^3 B_{2,3}^{-1}a^{-3} \rho_{3}= B_{2,3}\rho_{3}\label{eq:conjrho3}\\
\rho_{2} B_{1,2} \rho_{3}^{-1}\ldotp B_{2,3} \ldotp \rho_{3} B_{1,2}^{-1} \rho_{2}^{-1}&= b\garside[3]a^{-1} \ldotp B_{2,3} \ldotp a \garside[3]^{-1}b^{-1}= b\garside[3](a^{-1} B_{2,3} a) \garside[3]^{-1}b^{-1} \notag\\
&=b a \sigma_{1}^2 a^{-1} b^{-1}=\sigma_{2} a^2 \sigma_{1}^2 a^{-2} \sigma_{2}^{-1}\quad \text{from \req{defab}}\ \notag\\
&= B_{2,3}^{-1} \quad \text{from \req{cyclicperm}} .\label{eq:conjrho3a}
\end{align}
Here we have used \req{conjgarside}, as well as the standard property of $\garside$ (in $B_{n}$) that $\garside \sigma_{i} \garside^{-1}=\sigma_{n-i}$ for all $1\leq i\leq n-1$. So under conjugation by $b\garside[3]a^{-1}$,
\begin{equation}\label{eq:actconj2}
\left\{\begin{aligned}
\rho_{3} & \mapsto B_{2,3}\rho_{3}\\
B_{2,3}\rho_{3} & \mapsto \rho_{3}.
\end{aligned}\right.
\end{equation}
Relations~\reqref{actconj1} and~\reqref{actconj2} thus describe the action of $\quat$ on $\F[2]\left(\rho_{3},B_{2,3}\right)$, from which we may easily deduce its action on $\F[3]\left(\rho_{3}^2,B_{2,3}, \rho_{3}B_{2,3}\rho_{3}^{-1}\right)$:
\begin{align*}
\rho_{2}\rho_{1} \ldotp \rho_{3}^2 \ldotp \rho_{1}^{-1}\rho_{2}^{-1} &= \rho_{3}^{-2}\\
\rho_{2}\rho_{1} \ldotp B_{2,3} \ldotp \rho_{1}^{-1}\rho_{2}^{-1} &= \rho_{3}^{-2}\ldotp \rho_{3} B_{2,3}^{-1}\rho_{3}^{-1} \ldotp \rho_{3}^2\\
\rho_{2}\rho_{1} \ldotp \rho_{3}B_{2,3}\rho_{3}^{-1} \ldotp \rho_{1}^{-1}\rho_{2}^{-1} &=
\rho_{3}^{-2}\ldotp B_{2,3}^{-1} \ldotp \rho_{3}^2\\
\rho_{2} B_{1,2} \rho_{3}^{-1}\ldotp \rho_{3}^2 \ldotp \rho_{3} B_{1,2}^{-1} \rho_{2}^{-1}&= B_{2,3}\ldotp \rho_{3}B_{2,3}\rho_{3}^{-1} \ldotp \rho_{3}^2 \\
\rho_{2} B_{1,2} \rho_{3}^{-1}\ldotp B_{2,3} \ldotp \rho_{3} B_{1,2}^{-1} \rho_{2}^{-1}&=  B_{2,3}^{-1}\\
\rho_{2} B_{1,2} \rho_{3}^{-1}\ldotp \rho_{3}B_{2,3}\rho_{3}^{-1} \ldotp \rho_{3} B_{1,2}^{-1} \rho_{2}^{-1}&= B_{2,3} \ldotp \rho_{3} B_{2,3}^{-1} \rho_{3}^{-1} \ldotp B_{2,3}^{-1}
\end{align*}

We now record the action of $\Z_{3}$ on $\F[3]\rtimes \quat$. Since
\begin{equation*}
\rho_{3}^{-2}=\sigma_{2}\sigma_{1}^2\sigma_{2}=  \sigma_{2}\sigma_{1}^2\sigma_{2}^{-1} \sigma_{2}^{-2}=B_{1,3}B_{2,3}
\end{equation*}
and
\begin{equation}\label{eq:rho2rho1}
(\rho_{2}\rho_{1})^2=(\rho_{2} B_{1,2} \rho_{3}^{-1})^2=b^4=\garside[3]^2=B_{1,2}B_{1,3}B_{2,3},
\end{equation}
we have
\begin{equation}\label{eq:exprb12}
B_{1,2}= (\rho_{2}\rho_{1})^2 B_{2,3}^{-1} B_{1,3}^{-1}= (\rho_{2}\rho_{1})^2 \rho_{3}^2.
\end{equation}
So
\begin{align*}
a^4 \ldotp \rho_{2}\rho_{1} \ldotp a^{-4} &= \rho_{1}^{-1}\rho_{3} \quad \text{by \req{cyclicperm}}\\
&= \rho_{1}^{-1}\rho_{2}^{-1}\ldotp \rho_{2} B_{1,2} \rho_{3}^{-1} \ldotp \rho_{3}^2 B_{1,2}^{-1}\\
& =(\rho_{2}\rho_{1})^{-1} \ldotp \rho_{2} B_{1,2} \rho_{3}^{-1}\ldotp  (\rho_{2}\rho_{1})^{-2} \quad\text{by \req{exprb12}}\\
&= \rho_{2}\rho_{1} \ldotp \rho_{2} B_{1,2} \rho_{3}^{-1} \quad \text{by \req{rho2rho1}},
\end{align*}
and using the fact that $(\rho_{2}\rho_{1})^2=\ft[3]$ is of order $2$ and is central in $B_{3}(\rp)$. Further,
\begin{align*}
a^4 \ldotp \rho_{2} B_{1,2} \rho_{3}^{-1} \ldotp a^{-4} &=a^4 \ldotp \rho_{2}  \rho_{3}^{-1} B_{1,2}\ldotp a^{-4}= \rho_{1}^{-1} \rho_{2} a^{-1} \sigma_{2}^2 a \quad\text{by \req{cyclicperm}}\\
&= \rho_{1}^{-1} \rho_{2} \rho_{1}^{-1} \sigma_{1} \sigma_{2} \sigma_{2}^2 \sigma_{2}^{-1} \sigma_{1}^{-1} \rho_{1} \quad \text{by \req{defab}}\\
&= \rho_{1}^{-1} \rho_{2} \rho_{1}^{-1} \sigma_{1} \sigma_{2}^2  \sigma_{1}\sigma_{1}^{-2} \rho_{1}= \rho_{1}^{-1} \ldotp \rho_{2} \rho_{1} \sigma_{1}^{-2} \ldotp \rho_{1}= \rho_{1}^{-1} \rho_{1} \rho_{2} \rho_{1}=\rho_{2} \rho_{1}
\end{align*}
using \repr{present}. This describes the action of $\Z_{3}$ on the $\quat$-factor. As for the action of $\Z_{3}$ on $\F[2](B_{2,3},\rho_{3})$, using relations~\reqref{rho2rho1} and~\reqref{exprb12}, we have
\begin{align*}
a^4 B_{2,3} a^{-4}&= B_{1,2}^{-1}\quad\text{by \req{cyclicperm}}\\
&= \ft[3] \rho_{3}^{-2}\quad\text{by equations~\reqref{rho2rho1} and \reqref{exprb12},}
\end{align*}
and
\begin{align*}
a^4 \rho_{3} a^{-4}&=\rho_{2}^{-1}\quad\text{by \req{cyclicperm}}\\
&= \rho_{3}^{-1} B_{1,2} \ldotp B_{1,2}^{-1} \rho_{3} \rho_{2}^{-1}= \rho_{3}^{-1}(\rho_{2}\rho_{1})^2 \rho_{3}^2 (\rho_{2} B_{1,2} \rho_{3}^{-1})^{-1} \quad\text{by \req{exprb12}}\\
&= \rho_{3} (\rho_{2} B_{1,2} \rho_{3}^{-1}) \quad\text{by \req{rho2rho1}}
\end{align*}
Hence the action on $\F[3]\left(\rho_{3}^2,B_{2,3}, \rho_{3}B_{2,3}\rho_{3}^{-1}\right)$ is given by
\begin{align*}
a^4 \rho_{3}^2 a^{-4}=& \rho_{3} (\rho_{2} B_{1,2} \rho_{3}^{-1}) \rho_{3} (\rho_{2} B_{1,2} \rho_{3}^{-1})\\
&=\rho_{3} (\rho_{2} B_{1,2} \rho_{3}^{-1}) \rho_{3} (\rho_{2} B_{1,2} \rho_{3}^{-1})^{-1}\ft[3] \quad\text{by \req{rho2rho1}}\\
&=\rho_{3}  B_{2,3}\rho_{3} \ft[3]\quad\text{by \req{conjrho3}}\\
&= \ft[3] \ldotp \rho_{3} B_{2,3} \rho_{3}^{-1}\ldotp \rho_{3}^{2}\\
a^4 B_{2,3} a^{-4}&= \garside[3]^2 \ldotp  \rho_{3}^{-2}\\
a^4 \rho_{3}B_{2,3} \rho_{3}^{-1} a^{-4}&=
\rho_{3} (\rho_{2} B_{1,2} \rho_{3}^{-1})\ft[3] \rho_{3}^{-2} (\rho_{2} B_{1,2} \rho_{3}^{-1})^{-1}\rho_{3}^{-1} \\
&= \ft[3] \rho_{3} (B_{2,3}\rho_{3})^{-2}\rho_{3}^{-1} \quad\text{by \req{conjrho3}}\\
&= \ft[3]  B_{2,3}^{-1}\ldotp \rho_{3}^{-2}\ldotp \rho_{3} B_{2,3}^{-1}\rho_{3}^{-1}
\end{align*}
Setting $u=a^4$, $x=\rho_{2}\rho_{1}$, $y=\rho_{2} B_{1,2} \rho_{3}^{-1}$, $z_{1}=\rho_{3}^2$, $z_{2}=B_{2,3}$ and $z_{3}=\rho_{3}B_{2,3}\rho_{3}^{-1}$ yields the desired actions, and completes the proof of \repr{gamma2rp34}.
\end{proof}

\section{The derived series of $B_n(\rp)$}\label{sec:dsrp2}

In this section, we study the derived series of $B_n(\rp)$ and prove Theorems~\ref{th:dsbn} and~\ref{th:dsb4}. We start by showing that for all $n\neq 2,3,4$, $(B_n(\rp))^{(1)}$ is perfect. We then study the cases $n=3,4$ in more detail. If $n=3$, we are able to determine completely the derived series of $B_3(\rp)$, and deduce that it is residually soluble. If $n=4$, in \reth{dsb4} we obtain some partial results on the derived series of $B_4(\rp)$ and its quotients.

\begin{proof}[Proof of \reth{dsbn}]
Cases~(\ref{it:ds1}) and~(\ref{it:ds2}) follow directly from the fact that
$B_1(\rp)\cong \Z_{2}$ and $B_2(\rp)\cong \quat[16]$.
Now consider case~(\ref{it:ds4}), i.e.\ $n\geq 5$. 
Let $H\subseteq (B_n(\rp))^{(1)}$ be a normal subgroup of $B_n(\rp)$ such that
$A=(B_n(\rp))^{(1)}/H$ is Abelian (notice that this condition is satisfied if 
$H=(B_n(\rp))^{(2)}$), and let
\begin{equation*}
\left\{ \begin{aligned}
\pi \colon\thinspace B_n(\rp) & \to B_n(\rp)/H\\
\beta &\mapsto \widetilde{\beta}
\end{aligned}\right.
\end{equation*}
denote the canonical projection. Then the Abelianisation homomorphism $\alpha$ of \req{ses} factors through $B_n(\rp)/H$, in other words there exists a (surjective) homomorphism $\map{\widehat{\alpha}}{B_n(\rp)/H}[\bnab{n}]$ satisfying $\alpha= \widehat{\alpha}\circ \pi$. So \req{ses} induces the following short exact sequence:
\begin{equation*}
1\to A \to B_n(\rp)/H \stackrel{\widehat{\alpha}}{\to} \bnab{n} \to 1.
\end{equation*} 
In particular, $\bnab{n}\cong \Z_{2}\oplus \Z_{2}$ is a quotient of $B_n(\rp)/H$. We claim that the two are in fact isomorphic, which using the above short exact sequence will imply that $B_n(\rp))^{(1)}=H$, and thus $(B_n(\rp))^{(1)}=(B_n(\rp))^{(2)}$. To prove the claim, first note that $\widehat{\sigma_1},\ldots, \widehat{\sigma_{n-1}}, \widehat{\rho_1},\ldots, \widehat{\rho_n}$ generate $B_n(\rp)/H$. Since $\alpha(\sigma_i)=\alpha(\sigma_1)$ for all $1\leq i\leq n-1$, it follows that
$\widehat{\alpha}(\widehat{\sigma_i})= \widehat{\alpha}(\widehat{\sigma_1})$. So there exist
$t_i\in A$, with $t_{1}=1$, such that $\widehat{\sigma_i}=t_i\widehat{\sigma_1}$.

We now apply $\pi$ to each of the relations of \repr{present}. First suppose that
$3\leq i\leq n-1$. Since $\sigma_i$ commutes with $\sigma_1$, we have that
$\widehat{\sigma_1}\cdot t_i\widehat{\sigma_1} =t_i\widehat{\sigma_1}\cdot \widehat{\sigma_1}$,
and hence $t_i$ commutes with $\widehat{\sigma_1}$.

Now let $4\leq i\leq n-1$ (such an $i$ exists since $n\geq 5$). Since $\sigma_i$ commutes with $\sigma_2$, we obtain $t_i\widehat{\sigma_1}\cdot t_2\widehat{\sigma_1} =t_2\widehat{\sigma_1}\cdot
t_i\widehat{\sigma_1}$. But $A$ is Abelian, and so it follows from the previous paragraph that $t_2$ commutes with $\widehat{\sigma_1}$. Applying this to the image under $\pi$ of the relation $\sigma_1\sigma_2\sigma_1=\sigma_2\sigma_1\sigma_2$, we see that $t_2=t_2^2$, and hence $t_2=1$.

Next, if $i\geq 2$ then the relation $\sigma_i\sigma_{i+1}\sigma_i=
\sigma_{i+1}\sigma_i\sigma_{i+1}$ implies that $t_i=t_{i+1}$, and so
$t_2=\dotsc=t_{n-1}=1$. Hence $\widehat{\sigma_1}=\widehat{\sigma_2}=\dotsc
=\widehat{\sigma_{n-1}}$, and we denote this common element by $\sigma$.

Let $1\leq j \leq n$. Then from \repr{present} there exists $1\leq i \leq n-1$ such that $\rho_j$ and $\sigma_i$ commute. So in the quotient $B_n(\rp)/H$, $\widehat{\rho_j}$ commutes with $\sigma$ for all $1\leq j\leq n$. If $1\leq i\leq n-1$, the relation $\rho_{i+1}=\sigma_i^{-1}\rho_{i}\sigma_i^{-1}$ implies that
$\widehat{\rho_{i+1}}= \widehat{\rho_i}\sigma^{-2}$. Hence $\widehat{\rho_{i+1}}=\widehat{\rho_{1}}\sigma^{-2i}$. Projecting the relations $\rho_{i+1}^{-1}\rho_{i}^{-1}\rho_{i+1}\rho_{i}= \sigma_i^2$ into $B_n(\rp)/H$, where $1\leq i\leq n-1$, we obtain $\sigma^2=1$, and so
$\widehat{\rho_i}=\widehat{\rho_1}$ for all $1 \leq i\leq n$. Finally, by projecting the surface relation of $B_{n}(\rp)$ into $B_n(\rp)/H$, $\widehat{\rho_{1}}^2=\sigma^{2(n-1)}=1$. Therefore the group  $B_n(\rp)/H$ is a quotient of $\Z_2 \oplus \Z_2$. But we know already that $\bnab{n}\cong \Z_{2}\oplus \Z_{2}$ is a quotient of $B_n(\rp)/H$, and so this proves the claim. Taking $H=(B_n(\rp))^{(2)}$, it follows that the group  $(B_n(\rp))^{(1)}$ is perfect. A presentation of $(B_n(\rp))^{(1)}$ will be given in \repr{fullpres}. This proves part~(\ref{it:ds4}).

We now consider case~(\ref{it:ds3}), so $n=3$. Parts~(\ref{it:ds3gam}) and~(\ref{it:ds3a}) are just restatements of the results obtained in \repr{gamma2rp34}.
To prove part~(\ref{it:ds3b}), one may check easily using the presentation of \repr{gamma2rp34} that the Abelianisation $(B_{3}(\rp))^{(1)}/(B_{3}(\rp))^{(2)}$ of $(B_{3}(\rp))^{(1)}$ is cyclic of order $3$, generated by the Abelianisation of $a^4$. Since $(B_{3}(\rp))^{(1)}$ is isomorphic to $(\F[3]\rtimes \quat)\rtimes \Z_{3}$, where the $\Z_{3}$-factor is generated by $a^4$, we obtain $(B_{3}(\rp))^{(2)}\cong \F[3]\rtimes \quat$, where the action is once more given by \repr{gamma2rp34}. To see that the quotient $B_{3}(\rp)/(B_{3}(\rp))^{(2)}$ is isomorphic to $\dih{12}$, note first that we have the following commutative diagram of short exact sequences:
\begin{equation*}
\xymatrix{%
1 \ar[r] & (B_{3}(\rp))^{(1)} \ar[d]\ar[r] & B_{3}(\rp) \ar[d]\ar[r] & \gpab[(B_{3}(\rp))] \ar@{=}[d]\ar[r] & 1 \\
1 \ar[r] & (B_{3}(\rp))^{(1)}/(B_{3}(\rp))^{(2)} \ar[r] & B_{3}(\rp)/(B_{3}(\rp))^{(2)} \ar[r] & \gpab[(B_{3}(\rp))]\ar[r] & 1 }
\end{equation*}
Since $(B_{3}(\rp))^{(1)}/(B_{3}(\rp))^{(2)}\cong \Z_{3}$ and $\gpab[(B_{3}(\rp))]\cong \Z_{2} \oplus \Z_{2}$, it follows that the quotient $B_{3}(\rp)/(B_{3}(\rp))^{(2)}$ is an extension of $\Z_{3}$ by $\Z_{2} \oplus \Z_{2}$. We claim that the action is non trivial. To see this, we consider the conjugate of $a^4$ (which is a coset representative of the generator of $(B_{3}(\rp))^{(1)}/(B_{3}(\rp))^{(2)}$) by $\sigma_{1}$ (which is a coset representative of $\overline{\sigma}\in \bnab{3}$):
\begin{align*}
\sigma_{1}a^4\sigma_{1}^{-1}&= \sigma_{1}(a^4\sigma_{1}^{-1} a^{-4}) a^4
= \sigma_{1}a^{-1}\sigma_{2}^{-1} a a^4\quad \text{by \req{cyclicperm}}\\
&=\sigma_{1}\rho_{1}^{-1} \sigma_{1}\sigma_{2}\sigma_{2}^{-1} a^5\quad \text{by \req{defab}}\\
&=\rho_{2}^{-1} \rho_{1}^{-1}\rho_{2}^{-1}\rho_{3}^{-1}a^8\quad \text{by \repr{present} and \req{a3b2}}.
\end{align*}
Now 
\begin{align*}
\rho_{2}^{-1} \rho_{1}^{-1}\rho_{2}^{-1}\rho_{3}^{-1}&=\rho_{3}(\rho_{3}^{-1}\rho_{2}^{-1} \rho_{1}^{-1}\rho_{2}^{-1})\rho_{3}^{-1}=\rho_{3}\left(\rho_{3}^{-2}B_{1,2}
(\rho_{2}B_{1,2}\rho_{3}^{-1})^{-1} (\rho_{2}\rho_{1})^{-1}\right)\rho_{3}^{-1}\\
&= \rho_{3}\left((\rho_{2}\rho_{1})^2 
(\rho_{2}B_{1,2}\rho_{3}^{-1})^{-1} (\rho_{2}\rho_{1})^{-1}\right)\rho_{3}^{-1}\quad \text{by \req{exprb12}}\\
&= \rho_{3}\left( (\rho_{2}B_{1,2}\rho_{3}^{-1})^{-1} (\rho_{2}\rho_{1})\right)\rho_{3}^{-1} \quad \text{by \req{rho2rho1}}.
\end{align*}
But $(\rho_{2}B_{1,2}\rho_{3}^{-1})^{-1} (\rho_{2}\rho_{1}) \in (B_{3}(\rp))^{(2)}$from \repr{gamma2rp34}, and since $(B_{3}(\rp))^{(2)} \triangleleft B_{3}(\rp)$, it follows that $\rho_{2}^{-1} \rho_{1}^{-1}\rho_{2}^{-1}\rho_{3}^{-1}\in (B_{3}(\rp))^{(2)}$. Thus $\sigma_{1}a^4\sigma_{1}^{-1}$ is congruent modulo $(B_{3}(\rp))^{(2)}$ to $a^{-4}$, and the action of $\overline{\sigma}$ on $(B_{3}(\rp))^{(1)}/(B_{3}(\rp))^{(2)}$ is multiplication by $-1$. In particular, $B_{3}(\rp)/(B_{3}(\rp))^{(2)}$ is a non Abelian group of order $12$. Of the three non-Abelian groups of order $12$, $B_{3}(\rp)/(B_{3}(\rp))^{(2)}$ cannot be isomorphic to $\an[4]$ since the latter has no normal subgroup of order $3$. It cannot be isomorphic to $\dic{12}=\Z_{3} \rtimes \Z_{4}$ (with non-trivial action) either, since $\dic{12}$ has a unique subgroup of order $3$ with quotient $\Z_{4}$. We conclude that $B_{3}(\rp)/(B_{3}(\rp))^{(2)}\cong \dih{12}$. By the short exact sequence
\begin{multline*}
1 \to (B_{3}(\rp))^{(2)}/(B_{3}(\rp))^{(3)} \to B_{3}(\rp)/(B_{3}(\rp))^{(3)} \to\\
B_{3}(\rp)/(B_{3}(\rp))^{(2)}\to 1,
\end{multline*}
it follows that $B_{3}(\rp)/(B_{3}(\rp))^{(3)}$ is an extension of $\Z_{2}^4$ by $\dih{12}$, so is of order $192$. This proves part~(\ref{it:ds3b}).

Now let us prove~(\ref{it:ds3c}). The first part, that $(B_{3}(\rp))^{(2)}/(B_{3}(\rp))^{(3)}\cong \Z_{2}^4$, follows easily by Abelianising the presentation of $(B_{3}(\rp))^{(2)}\cong \F[3]\rtimes \quat$ given in \repr{gamma2rp34}. Letting $\phi$ denote the Abelianisation homomorphism, one observes that $\phi(z_{1})=(\overline{1}, \overline{0}, \overline{0}, \overline{0})$, $\phi(z_{2})= \phi(z_{3})= (\overline{0}, \overline{1}, \overline{0}, \overline{0})$, $\phi(x)= (\overline{0}, \overline{0}, \overline{1}, \overline{0})$, and $\phi(y)= (\overline{0}, \overline{0}, \overline{0}, \overline{1})$. The restriction of $\phi$ to $\F[3]$ is surjective onto the subgroup $H$ of $(B_{3}(\rp))^{(2)}/(B_{3}(\rp))^{(3)}\cong \Z_{2}^4$ generated by the $\phi(z_{i})$, so $H\cong \Z_{2}^2$. The quotient $Q$ of $(B_{3}(\rp))^{(2)}/(B_{3}(\rp))^{(3)}$ by $H$ is thus isomorphic to the subgroup of $(B_{3}(\rp))^{(2)}/(B_{3}(\rp))^{(3)}$ generated by $\phi(x)$ and $\phi(y)$, and so is also isomorphic to $\Z_{2}^2$. Since $\quat=\ang{x,y}$, $\phi$ induces a surjective homomorphism $\map{\overline{\phi}}{\quat}[Q]$ whose kernel is $\ang{x^2}$. But as an element of $(B_{3}(\rp))^{(2)}$, $x^2=\ft[3]\in \ker{\phi}$, and denoting $\ker{\phi\left\lvert_{\F[3]}\right.}$ by $L$, we obtain the following commutative diagram of short exact sequences:
\begin{equation*}
\xymatrix{%
& 1\ar[d] & 1\ar[d] & 1\ar[d] & \\
1 \ar[r] & \ar[r] L \ar[d] & (B_{3}(\rp))^{(3)}\ar[d] \ar[r]  & \ang{\ft[3]} \ar[d]\ar[r] & 1\\
1 \ar[r] & \ar[r] \F[3] \ar[d]^{\phi\left\lvert_{\F[3]}\right.} & (B_{3}(\rp))^{(2)}\ar[r] \ar[d]^{\phi} & \quat\ar[d] \ar[r] & 1\\
1 \ar[r] & \ar[r] H \ar[d] & (B_{3}(\rp))^{(2)}/(B_{3}(\rp))^{(3)} \ar[d] \ar[r]  & Q \ar[d]\ar[r] & 1.\\
& 1 & 1 & 1 & }
\end{equation*}
Since $\ft[3]=x^2\in (B_{3}(\rp))^{(3)}$, it follows that the upper short exact sequence splits, and the fact that $\ang{\ft[3]}$ is central implies that the splitting gives rise to a direct product. We conclude that $(B_{3}(\rp))^{(3)}\cong L\oplus \Z_{2}$. Now $L$ is the kernel of the homomorphism $\map{\phi\left\lvert_{\F[3]}\right.}{\F[3](z_{1},z_{2},z_{3})}[H]$ which under identification of $H$ with $\Z_{2}\oplus \Z_{2}$ sends $z_{1}$ to $(\overline{1}, \overline{0})$, say, and $z_{2}$ and $z_{3}$ to $(\overline{0}, \overline{1})$. An application of the Reidemeister-Schreier rewriting process shows that $L$ is a free group $\F[9]$ of rank $9$. Thus $(B_{3}(\rp))^{(3)}\cong \F[9]\oplus \Z_{2}$ and $(B_{3}(\rp))^{(3)}/(B_{3}(\rp))^{(4)}\cong \Z^9\oplus \Z_{2}$. It is then clear that $(B_{3}(\rp))^{(i)}\cong (\F[9])^{(i-3)}$ for all $i\geq 4$. From the short exact sequence
\begin{multline*}
1 \to (B_{3}(\rp))^{(3)}/(B_{3}(\rp))^{(4)} \to B_{3}(\rp)/(B_{3}(\rp))^{(4)} \to\\ 
B_{3}(\rp)/(B_{3}(\rp))^{(3)} \to 1,
\end{multline*}
we see that $B_{3}(\rp)/(B_{3}(\rp))^{(4)}$ is an extension of $\Z^9\oplus \Z_{2}$ by a group of order $192$, so is infinite. This proves part~(\ref{it:ds3c}), and completes the proof of \reth{dsbn}.
\end{proof}

We obtain easily the following corollary of \reth{dsbn}:
\begin{cor}\label{cor:residb4}
Let $n\in\N$, $n\neq 4$. Then $B_n(\rp)$ is residually soluble if and only if $n\leq 3$.\hfill\qedsymbol
\end{cor}

We now turn our attention to the remaining case, $n=4$.

\begin{proof}[Proof of \reth{dsb4}]
Part~(\ref{it:ds4gam}) follows from the first paragraph of the proof of \repr{gamma2rp34}. So let us prove part~(\ref{it:ds4a}). For this, we shall study the following presentation of the group $(B_4(\rp))^{(1)}$ which may be deduced from \repr{fullpres} (the notation $\alpha,\beta$ etc.\ is that of \repr{fullpres}):
\begin{enumerate}
\item[\underline{\textbf{generators:}}]\mbox{}

\vspace{-0.7cm}

\noindent
\begin{minipage}[t]{0.35\linewidth}
\begin{align*}
& B_{1}=\eta_{2}= \rho_{2}\sigma_{1}^{-1}\rho_{1}^{-1}\sigma_{1}^{-1}\\
& B_{2}=\kappa_{2}= \sigma_{1}\rho_{2}\rho_{1}^{-1}\sigma_{1}^{-1}\\
& B_{3}=\theta_{2}= \sigma_{1}\rho_{1} \rho_{2}\sigma_{1}^{-1}\\
& B_{4}=\lambda_{2}= \sigma_{1}\rho_{1} \sigma_{1}\rho_{2}\\
& Y_{1}=\alpha_{2}= \sigma_{2}\sigma_{1}^{-1}\\ 
& Y_{2}=\beta_{2}= \sigma_{1}\sigma_{2}\\
& Y_{3}=\gamma_{2}= \sigma_{1}\rho_{1} \sigma_{2}\sigma_{1}^{-1} \rho_{1}^{-1}\sigma_{1}^{-1}\\
& Y_{4}=\tau_{2}= \sigma_{1}\rho_{1} \sigma_{1}\sigma_{2} \rho_{1}^{-1}\sigma_{1}^{-1} 
\end{align*}
\end{minipage}
\begin{minipage}[t]{0.35\linewidth}
\begin{align*}
& C_{1}=\eta_{3}= \rho_{3}\sigma_{1}^{-1}\rho_{1}^{-1}\sigma_{1}^{-1}\\
& C_{2}=\kappa_{3}= \sigma_{1}\rho_{3}\rho_{1}^{-1}\sigma_{1}^{-1}\\
& C_{3}=\theta_{3}= \sigma_{1}\rho_{1} \rho_{3}\sigma_{1}^{-1} \\  
& C_{4}=\lambda_{3}= \sigma_{1}\rho_{1} \sigma_{1}\rho_{3}\\
& Z_{1}=\alpha_{3}= \sigma_{3}\sigma_{1}^{-1}\\
& Z_{2}=\beta_{3}= \sigma_{1}\sigma_{3}\\
& Z_{3}=\gamma_{3}= \sigma_{1}\rho_{1} \sigma_{3}\sigma_{1}^{-1} \rho_{1}^{-1}\sigma_{1}^{-1}\\
& Z_{4}=\tau_{3}= \sigma_{1}\rho_{1} \sigma_{1}\sigma_{3} \rho_{1}^{-1}\sigma_{1}^{-1}
\end{align*}
\end{minipage}
\begin{minipage}[t]{0.3\linewidth}
\begin{align*}
& D_{1}=\eta_{4}= \rho_{4}\sigma_{1}^{-1}\rho_{1}^{-1}\sigma_{1}^{-1} \\
& D_{2}=\kappa_{4}= \sigma_{1}\rho_{4}\rho_{1}^{-1}\sigma_{1}^{-1}\\ 
& D_{3}=\theta_{4}= \sigma_{1}\rho_{1} \rho_{4}\sigma_{1}^{-1} \\
& D_{4}=\lambda_{4}= \sigma_{1}\rho_{1} \sigma_{1}\rho_{4}\\
& A_{1}=\eta_{1}= \rho_{1}\sigma_{1}^{-1}\rho_{1}^{-1}\sigma_{1}^{-1}\\  
& A_{3}=\theta_{1}= \sigma_{1}\rho_{1} \rho_{1}\sigma_{1}^{-1}\\
& A_{4}=\lambda_{1}= \sigma_{1}\rho_{1} \sigma_{1}\rho_{1}\\
& X_{2}=\beta_{1}= \sigma_{1}^2\\
& X_{4}=\tau_{1}= \sigma_{1}\rho_{1} \sigma_{1}^2 \rho_{1}^{-1}\sigma_{1}^{-1}.
\end{align*}
\end{minipage}
For notational reasons, we set $A_{2}=X_{1}=X_{3}=1$.

\medskip

\item[\underline{\textbf{relators:}}]\mbox{}

\noindent
\begin{minipage}{0.33\linewidth}
\begin{align}
& Z_{2} X_{2}^{-1}Z_{1}^{-1}\label{eq:firstrel}\\
& X_{2} Z_{1} Z_{2}^{-1}\\
& Z_{4} X_{4}^{-1}Z_{3}^{-1} \label{eq:firstrelx4}\\
& X_{4} Z_{3} Z_{4}^{-1}\\
& Y_{2}Y_{1}^{-1}X_{2}^{-1}Y_{1}^{-1} \label{eq:y2y1sq}\\ 
& X_{2} Y_{1} X_{2} Y_{2}^{-2}  \label{eq:y1y2sq}\\ 
& Y_{4}Y_{3}^{-1}X_{4}^{-1}Y_{3}^{-1} \label{eq:y4y3sq}
\end{align}
\end{minipage}
\begin{minipage}{0.33\linewidth}
\begin{align}
& X_{4} Y_{3} X_{4} Y_{4}^{-2}\label{eq:y3y4sq}\\
& Y_{1}Z_{2}Y_{1}Z_{1}^{-1}Y_{2}^{-1}Z_{1}^{-1} \label{eq:y1z2}\\
& Y_{2} Z_{1} Y_{2} Z_{2}^{-1} Y_{1}^{-1} Z_{2}^{-1} \label{eq:y2z1}\\
& Y_{3}Z_{4}Y_{3}Z_{3}^{-1}Y_{4}^{-1}Z_{3}^{-1} \\
& Y_{4} Z_{3} Y_{4} Z_{4}^{-1} Y_{3}^{-1} Z_{4}^{-1}\\
& C_{2}X_{4}^{-1}C_{1}^{-1} \label{eq:x4c1c2}\\ 
& X_{2} C_{1} C_{2}^{-1}\label{eq:x2c1c2}
\end{align}
\end{minipage}
\noindent
\begin{minipage}{0.33\linewidth}
\begin{align}
& C_{4} X_{2}^{-1} C_{3}^{-1}\label{eq:x2c3c4}\\
& X_{4} C_{3}  C_{4}^{-1} \label{eq:x4c3c4}\\
& D_{2}X_{4}^{-1}D_{1}^{-1} \label{eq:x4d1d2}\\
& X_{2} D_{1} D_{2}^{-1} \label{eq:x2d1d2}\\
& D_{4} X_{2}^{-1} D_{3}^{-1} \label{eq:x2d3d4}\\
& X_{4} D_{3}  D_{4}^{-1}\label{eq:x4d3d4}\\
& Y_{1}Y_{4}^{-1}A_{1}^{-1} \label{eq:a1y1y4}
\end{align}
\end{minipage}

\noindent
\begin{minipage}{0.33\linewidth}
\begin{align}
& Y_{2} A_{1} Y_{3}^{-1}  \label{eq:a1y2y3}\\
& Y_{3} A_{4} Y_{2}^{-1} A_{3}^{-1} \\
& Y_{4} A_{3} Y_{1}^{-1} A_{4}^{-1} \\
& Y_{1}D_{2}Y_{4}^{-1}D_{1}^{-1} \\
& Y_{2} D_{1} Y_{3}^{-1} D_{2}^{-1} \\
& Y_{3} D_{4} Y_{2}^{-1} D_{3}^{-1} \\
& Y_{4} D_{3} Y_{1}^{-1} D_{4}^{-1} 
\end{align}
\end{minipage}
\noindent
\begin{minipage}{0.33\linewidth}
\begin{align}
& Z_{1}Z_{4}^{-1}A_{1}^{-1} \label{eq:a1z1z4}\\
& Z_{2} A_{1} Z_{3}^{-1}  \label{eq:a1z2z3}\\
& Z_{3} A_{4} Z_{2}^{-1} A_{3}^{-1} \label{eq:a3z2z3}\\
& Z_{4} A_{3} Z_{1}^{-1} A_{4}^{-1} \label{eq:a4z1z4}\\
& Z_{1}B_{2}Z_{4}^{-1}B_{1}^{-1} \\
& Z_{2} B_{1} Z_{3}^{-1} B_{2}^{-1} \\
& Z_{3} B_{4} Z_{2}^{-1} B_{3}^{-1} 
\end{align}
\end{minipage}
\begin{minipage}{0.33\linewidth}
\begin{align}
& Z_{4} B_{3} Z_{1}^{-1} B_{4}^{-1}\\
& B_{1} X_{4} X_{2} \label{eq:b1x2x4}\\ 
& B_{2}A_{1}^{-1}\label{eq:a1b2}\\
& X_{4}^{-1} A_{4} X_{2}^{-1} B_{3}^{-1} \label{eq:x4inva4}\\
& A_{3} B_{4}^{-1} \label{eq:b4a3}\\ 
& Y_{2}^{-1} B_{2} Y_{4}^{-1} C_{1}^{-1}\label{eq:y2invb2} \\
& Y_{1}^{-1} B_{1} Y_{3}^{-1} C_{2}^{-1} 
\end{align}
\end{minipage}

\noindent
\begin{minipage}[t]{0.33\linewidth}
\begin{align}
& Y_{4}^{-1} B_{4} Y_{2}^{-1} C_{3}^{-1} \label{eq:y4invb4}\\
& Y_{3}^{-1} B_{3} Y_{1}^{-1} C_{4}^{-1} \\
& Z_{1}^{-1} C_{1} Z_{3}^{-1} D_{2}^{-1}  \\
& Z_{2}^{-1} C_{2} Z_{4}^{-1} D_{1}^{-1} \label{eq:z2invc2}\\
& Z_{3}^{-1} C_{3} Z_{1}^{-1} D_{4}^{-1} \\
& Z_{4}^{-1} C_{4} Z_{2}^{-1} D_{3}^{-1} \label{eq:z4invc4}\\
& B_{4}^{-1} A_{1}^{-1} B_{1} A_{4}X_{2}^{-1}\label{eq:b4a1b1}\\
& B_{2} A_{3}  X_{2}^{-1} B_{3}^{-1} 
\end{align}
\end{minipage}
\begin{minipage}[t]{0.33\linewidth}
\begin{align}
& B_{2}^{-1} A_{3}^{-1} B_{3}  X_{4}^{-1}  \\
& B_{1}^{-1} A_{4}^{-1} B_{4} A_{1}  X_{4}^{-1} \\
& C_{4}^{-1} B_{1}^{-1} C_{1} B_{4}Y_{2}^{-1} Y_{1}^{-1} \\
& C_{3}^{-1} B_{2}^{-1} C_{2} B_{3} Y_{1}^{-1} Y_{2}^{-1} \label{eq:c3invb2inv}\\
& C_{2}^{-1} B_{3}^{-1} C_{3} B_{2} Y_{4}^{-1} Y_{3}^{-1} \\
& C_{1}^{-1} B_{4}^{-1} C_{4} B_{1} Y_{3}^{-1} Y_{4}^{-1} \\
& D_{4}^{-1} C_{1}^{-1} D_{1} C_{4}Z_{2}^{-1} Z_{1}^{-1}\\
& D_{3}^{-1} C_{2}^{-1} D_{2} C_{3} Z_{1}^{-1} Z_{2}^{-1}
\end{align}
\end{minipage}
\begin{minipage}[t]{0.33\linewidth}
\begin{align}
& D_{2}^{-1} C_{3}^{-1} D_{3} C_{2} Z_{4}^{-1} Z_{3}^{-1} \\
& D_{1}^{-1} C_{4}^{-1} D_{4} C_{1} Z_{3}^{-1} Z_{4}^{-1}\\ 
& Y_{2}Z_{1}Z_{2}Y_{1}X_{2} A_{4}^{-1}A_{1}^{-1}\label{eq:surf1}\\
& X_{2}Y_{1}Z_{2}Z_{1}Y_{2}X_{2} A_{3}^{-1}\label{eq:surf2}\\
& Y_{4}Z_{3}Z_{4}Y_{3}X_{4}  A_{3}^{-1}\\
& X_{4}Y_{3}Z_{4}Z_{3}Y_{4} A_{1}^{-1} A_{4}^{-1}.\label{eq:lastrel}
\end{align}
\end{minipage}
\end{enumerate}

We now Abelianise this presentation to deduce that $(B_4(\rp))^{(1)}/(B_4(\rp))^{(2)}\cong \Z_{3}$. We could do this directly, but it will be convenient for what follows to carry out a partial Abelianisation first. Let $\Lambda$ denote the group obtained from the above presentation of $(B_4(\rp))^{(1)}$ by adding the relations that the following generators commute pairwise: $A_{i},B_{i},C_{i}, D_{i}, X_{i}, Z_{i}$ for $i=1,\ldots,4$ i.e.\ all of the generators of $(B_4(\rp))^{(1)}$ commute pairwise, with the exception of the $Y_{i}$. From equations~\reqref{firstrel}, \reqref{x2c1c2}, \reqref{x2c3c4}, \reqref{x2d1d2} and~\reqref{x2d3d4}, we have
\begin{equation}\label{eq:x2zcd}
X_{2}=Z_{2}Z_{1}^{-1}= C_{2}C_{1}^{-1}= D_{2}D_{1}^{-1}=C_{3}^{-1}C_{4}=D_{3}^{-1}D_{4}
\end{equation}
and from equations~\reqref{firstrelx4}, \reqref{x4c1c2},  \reqref{x4c3c4}, \reqref{x4d1d2} and~\reqref{x4d3d4}, we have
\begin{equation}\label{eq:x4zcd}
X_{4}=Z_{4}Z_{3}^{-1}= C_{2}C_{1}^{-1}= D_{2}D_{1}^{-1}=C_{3}^{-1}C_{4}=D_{3}^{-1}D_{4}.
\end{equation}
So
\begin{equation}\label{eq:x2x4}
X_{2}=X_{4}
\end{equation}
and
\begin{equation}\label{eq:z2z3z1z4}
Z_{2}Z_{3}=Z_{1}Z_{4}.
\end{equation}
Now from equations~\reqref{a1y1y4}, \reqref{a1y2y3}, \reqref{a1z1z4}, \reqref{a1z2z3} and \reqref{a1b2}, we have
\begin{equation}\label{eq:a1y1y4b}
A_{1}=B_{2}=Y_{1}Y_{4}^{-1}=Y_{2}^{-1}Y_{3}=Z_{1}Z_{4}^{-1}=Z_{2}^{-1}Z_{3},
\end{equation}
hence
\begin{equation}\label{eq:z1z2z3z4}
Z_{1}Z_{2}=Z_{3}Z_{4}.
\end{equation}
Multiplying equations~\reqref{z2z3z1z4} and~\reqref{z1z2z3z4} yields
\begin{equation}\label{eq:zi2}
Z_{1}^2=Z_{3}^2,\quad Z_{2}^2=Z_{4}^2.
\end{equation}
Further, from equations~\reqref{a3z2z3}, \reqref{a4z1z4}, \reqref{b1x2x4}, \reqref{b4a3}, \reqref{x2zcd}, \reqref{x4zcd} and~\reqref{z1z2z3z4}, we have
\begin{align}
B_{4}&=A_{3}=Z_{3}A_{4}Z_{2}^{-1}= Z_{1}Z_{4}^{-1}A_{4}\label{eq:a3b4}\\
B_{1}&=X_{2}^{-1}X_{4}^{-1}=Z_{1}Z_{2}^{-1} Z_{3}Z_{4}^{-1}=Z_{1}^2 Z_{4}^{-2}.\label{eq:b1x2}
\end{align}
Substituting equations~\reqref{a1y1y4b},~\reqref{a3b4} and~\reqref{b1x2} into \req{b4a1b1} yields:
\begin{equation*}
X_{2}=B_{4}^{-1} A_{1}^{-1} B_{1} A_{4}= A_{4}^{-1} Z_{4} Z_{1}^{-1} Z_{4} Z_{1}^{-1} Z_{1}^2 Z_{4}^{-2}A_{4}=1.
\end{equation*}
Hence
\begin{gather}
X_{2}=X_{4}=1,\; Z_{1}=Z_{2},\; Z_{3}=Z_{4},\; Z_{1}^2=Z_{2}^2=Z_{3}^2=Z_{4}^2 \label{eq:x2x4b}\\
C_{1}=C_{2},\; D_{1}=D_{2},\; C_{3}=C_{4},\; D_{3}=D_{4},\; B_{1}=1\label{eq:c1c2}
\end{gather}
by equations~\reqref{x2zcd}, \reqref{x4zcd}, \reqref{x2x4}, \reqref{zi2} and~\reqref{b1x2}.
From equations~\reqref{y2y1sq}, \reqref{y1y2sq}, \reqref{y4y3sq} and~\reqref{y3y4sq}, we have $Y_{2}=Y_{1}^2$, $Y_{1}=Y_{2}^2$, $Y_{4}=Y_{3}^2$, $Y_{3}=Y_{4}^2$, so
\begin{equation}\label{eq:y1y2}
Y_{2}=Y_{1}^{-1},\; Y_{4}=Y_{3}^{-1},\; \text{and $Y_{i}^3=1$ for all $i=1,\ldots,4$.}
\end{equation}
Using equations~\reqref{y1z2}, \reqref{y2z1} and~\reqref{x2x4b}, we see that
\begin{equation}\label{eq:zyzyi}
1=Y_{1}Z_{2}Y_{1}Z_{1}^{-1}Y_{2}^{-1} Z_{1}^{-1}= Y_{1}Z_{1}Y_{1}Z_{1}^{-1}Y_{1}Z_{1}^{-1}
\end{equation}
and 
\begin{equation}\label{eq:zyzyii}
1=Y_{2}Z_{1}Y_{2}Z_{2}^{-1}Y_{1}^{-1} Z_{2}^{-1}= Y_{1}^{-1}Z_{1}Y_{1}^{-1} Z_{1}^{-1}Y_{1}^{-1} Z_{1}^{-1}.
\end{equation}
Inverting \req{zyzyii} and conjugating by $Z_{1}^{-1}$, we obtain
\begin{equation}\label{eq:zyzyiii}
1= Y_{1}Z_{1}Y_{1} Z_{1}^{-1}Y_{1}Z_{1}.
\end{equation}
Comparing equations~\reqref{zyzyi} and~\reqref{zyzyiii} yields 
\begin{equation}\label{eq:z12}
Z_{1}^2=1.
\end{equation}
From this and equations~\reqref{a1y1y4b} and~\reqref{x2x4b}, it follows that
\begin{equation}\label{eq:a12}
A_{1}^2=1.
\end{equation}
By equations~\reqref{x4inva4}, \reqref{surf1}, \reqref{a1y1y4b}, \reqref{surf2}, \reqref{a3b4}, \reqref{x2x4b}, \reqref{y1y2}, \reqref{z12} and \reqref{a12}, we see that
\begin{equation}\label{eq:a1a4}
A_{1}=A_{4}=B_{2}=B_{3},\; A_{3}=B_{4}=1.
\end{equation}
Now 
\begin{equation}\label{eq:c1y1y3}
C_{1}=Y_{2}^{-1} B_{2} Y_{4}^{-1}=Y_{1}A_{1}Y_{3}=Y_{1}^{-1} Y_{3}^{-1}
\end{equation}
by equations~\reqref{a1y1y4}, \reqref{y2invb2}, \reqref{y1y2} and~\reqref{c1y1y3}, and
\begin{equation}\label{eq:c3y3y1}
C_{3}=Y_{4}^{-1} B_{4} Y_{2}^{-1}= Y_{3}Y_{1}
\end{equation}
by equations~\reqref{y4invb4}, \reqref{y1y2} and \reqref{a1a4}, hence
\begin{equation}\label{eq:cequal}
C_{1}=C_{2}=C_{3}^{-1}=C_{4}^{-1}
\end{equation}
by \req{c1c2}.
Similarly,
\begin{equation}\label{eq:d1y3y1}
D_{1}=Z_{2}^{-1} C_{2} Z_{4}^{-1}=Z_{1}^{-1} Y_{1}^{-1} Y_{3}^{-1} Z_{3}^{-1}
\end{equation}
by equations~\reqref{z2invc2}, \reqref{x2x4b}, \reqref{c1c2} and \reqref{c1y1y3}, and
\begin{equation}\label{eq:d3y1y3}
D_{3}=Z_{4}^{-1} C_{4} Z_{2}^{-1}= Z_{3}^{-1} Y_{3} Y_{1} Z_{1}^{-1}
\end{equation}
by equations~\reqref{z4invc4}, \reqref{x2x4b}, \reqref{c1c2} and \reqref{c3y3y1}, hence
\begin{equation}\label{eq:dequal}
D_{1}=D_{2}=D_{3}^{-1}=D_{4}^{-1}
\end{equation}
by equations~\reqref{x2x4b}, \reqref{c1c2} and~\reqref{z12}.
Using equations~\reqref{c3invb2inv}, \reqref{y1y2}, \reqref{a12}, \reqref{a1a4}, \reqref{c3y3y1} and \reqref{cequal}, we see that
\begin{equation*}
1=Y_{1}^{-1}Y_{2}^{-1}=B_{3}^{-1} C_{2}^{-1} B_{2} C_{3}= A_{1}C_{2}^{-1} A_{1} C_{3}=C_{3}^2= (Y_{3}Y_{1})^2.
\end{equation*}
Hence $Y_{3}Y_{1}=Y_{1}^{-1}Y_{3}^{-1}$, and so
\begin{equation}\label{eq:callequal}
C_{1}=C_{2}=C_{3}=C_{4},\; D_{1}=D_{2}=D_{3}=D_{4}
\end{equation}
by equations~\reqref{cequal}, \reqref{d1y3y1}, \reqref{d3y1y3} and \reqref{dequal}, as well as the fact that $C_{2},C_{4},Z_{2}$ and $Z_{4}$ commute pairwise. We deduce also from equations~\reqref{dequal} and~\reqref{callequal} that $D_{i}^2=1$ for all $i=1,\ldots, 4$. Let $C$ (resp.\ $D$) denote the common value of the $C_{i} $ (resp.\ $D_{i}$), and let $A=A_{1}=A_{4}=B_{2}=B_{3}$. Running through the relations~\reqref{firstrel}--\reqref{lastrel} one by one, we see that our group $\Lambda$ has generators $Y_{1},Y_{3}, Z_{1}, Z_{3}, A,C$ and $D$ with the following defining relations:
\begin{equation}\label{eq:presd1d3}
\left\{
\begin{aligned}
& Y_{1}^3=Y_{3}^3=(Y_{1}Z_{1})^3=(Y_{3}Z_{3})^3=Z_{1}^2=Z_{3}^2=A^2=C^2=D^2=1\\
& A=CD=Y_{1}Y_{3}=Z_{1}Z_{3},\; C=Y_{3}Y_{1},\; Y_{1}D Y_{3}D=1,\; (Y_{3}Y_{1})^2=1 \\
& \text{$A,C,D,Z_{1}$ and $Z_{3}$ commute pairwise.}
\end{aligned}\right.
\end{equation}
Notice that we may write $D=Y_{3} Y_{1}Z_{3} Z_{1}$, and so $Y_{1},Y_{3}, Z_{1}, Z_{3}$ generate $\Lambda$.

If we now Abelianise $(B_4(\rp))^{(1)}$ completely by adding the relations that the $Y_{i}$ commute pairwise with each of the generators of $\Lambda$, we see that $A=C=D=Z_{1}=Z_{3}=1$, $Y_{3}=Y_{1}^{-1}, Y_{1}^3=1$, and thus $(B_4(\rp))^{(1)}/(B_4(\rp))^{(2)}\cong \Z_{3}$. We underline the fact that under the complete Abelianisation of $(B_4(\rp))^{(1)}$, the generators $A_{i},B_{i},C_{i},D_{i},X_{i},Z_{i}$, $i=1,\ldots,4$, of $(B_4(\rp))^{(1)}$ are sent to the trivial element, and $Y_{1},Y_{2}^{-1},Y_{3}^{-1}$ and $Y_{4}$ are sent to the same generator of $(B_4(\rp))^{(1)}/(B_4(\rp))^{(2)}$. Taking $b=\rho_{3}\sigma_{2} \sigma_{1}=\sigma_{2}^{-1} \sigma_{1}^{-1}\rho_{1}\in B_{4}(\rp)$ which we know to be of order $12$, consider $b^4$. Since $b^3=\rho_{3}\rho_{2}\rho_{1}$ by \rerem{defab}(\ref{it:defab}), we have
\begin{equation*}
b^4=\rho_{3}\rho_{2}\rho_{1} \ldotp \rho_{3}\sigma_{2} \sigma_{1}=C_{1}B_{4}A_{1}C_{4}Y_{1}X_{2}\in (B_4(\rp))^{(1)}.
\end{equation*}
Under Abelianisation, $b^4$ is thus sent to the $(B_4(\rp))^{(2)}$-coset of $Y_{1}$ which is a generator of $(B_4(\rp))^{(1)}/(B_4(\rp))^{(2)}$. Since $b^4$ is of order $3$, it follows that the short exact sequence
\begin{equation*}
1\to (B_4(\rp))^{(2)} \to (B_4(\rp))^{(1)} \to (B_4(\rp))^{(1)}/(B_4(\rp))^{(2)}\to 1
\end{equation*}
splits, and hence
\begin{equation*}
(B_4(\rp))^{(1)}\cong (B_4(\rp))^{(2)} \rtimes \Z_{3},
\end{equation*}
where the action on $(B_4(\rp))^{(2)}$ is given by conjugation by $b^4$. This proves part~(\ref{it:ds4a})(\ref{it:ds4ai}). To prove part~(\ref{it:ds4a})(\ref{it:ds4aii}), consider the short exact sequence
\begin{multline*}
1\to (B_4(\rp))^{(1)}/(B_4(\rp))^{(2)} \to B_4(\rp)/(B_4(\rp))^{(2)} \to\\ 
B_4(\rp)/(B_4(\rp))^{(1)}\to 1.
\end{multline*}
As in part~(\ref{it:ds3})(\ref{it:ds3b}) of the proof of \reth{dsbn}, since $(B_4(\rp))^{(1)}/(B_4(\rp))^{(2)}\cong \Z_{3}$ and $B_4(\rp)/(B_4(\rp))^{(1)}\cong \Z_{2}\oplus \Z_{2}$, to prove that $B_4(\rp)/(B_4(\rp))^{(2)}\cong \dih{12}$, it suffices to show that the action of $B_4(\rp)/(B_4(\rp))^{(1)}$ on the kernel is non trivial. To achieve this, notice that the action by conjugation of $\sigma_{1}$ (which is a representative of the generator $\overline{\sigma}$ of $B_4(\rp)/(B_4(\rp))^{(1)}$) on $Y_{1}$ (which from above is a representative of a generator of $(B_4(\rp))^{(1)}/(B_4(\rp))^{(2)}$) is given by $\sigma_{1} Y_{1} \sigma_{1}^{-1}=\sigma_{1}\sigma_{2}\sigma_{1}^{-2}= Y_{2}X_{2}^{-1}$.
Now modulo $(B_4(\rp))^{(2)}$, $Y_{2}X_{2}^{-1}$ is congruent to $Y_{2}$, which in turn is congruent to $Y_{1}^{-1}$. The action of $B_4(\rp)/(B_4(\rp))^{(1)}$ on $(B_4(\rp))^{(1)}/(B_4(\rp))^{(2)}$ is thus non trivial, which proves that $B_4(\rp)/(B_4(\rp))^{(2)}\cong \dih{12}$, and completes the proof of part~(\ref{it:ds4a})(\ref{it:ds4aii}).

To prove part~(\ref{it:ds4a})(\ref{it:ds4aiii}), let $n=4$ in the commutative diagram~\reqref{diaggam2} of short exact sequences. Recall that in the lower sequence, $\bnab{4}\cong \Z_{2}\oplus \Z_{2}$ is generated by two elements $\overline{\sigma}$ and $\overline{\rho}$, $\ker{\overline{\tau}}=\ang{\overline{\rho}}$, and $\overline{\tau}(\overline{\sigma})$, which we also denote by $\overline{\sigma}$, is the generator of the quotient $\left.\bnab{4}\right/\ang{\overline{\rho}}$. From the discussion following \req{diaggam2}, $K$ is of index $2$ in $P_{4}(\rp)$. Furthermore, the homomorphism $\alpha'$ sends the generator $B_{i,j}$, $1\leq i<j\leq 4$ (resp.\ $\rho_{k}$, $1\leq k\leq 4$) to the trivial element of $\ang{\overline{\rho}}$ (resp.\ to $\overline{\rho}$). This diagram may be continued vertically by taking commutator subgroups successively; in this way, we obtain the following commutative diagram of short exact sequences:
\begin{equation}\label{eq:bigdig2}
\begin{xy}*!C\xybox{%
\xymatrix{%
1 \ar[r] & \ar[r] K'' \ar[d] & (B_{4}(\rp))^{(3)}\ar[d] \ar[r]  & 1 \ar[d]\ar[r] & 1\\
1 \ar[r] & \ar[r] K' \ar[d] & (B_{4}(\rp))^{(2)}\ar[d] \ar[r]  & \Z_{2}\oplus \Z_{2} \ar[d]\ar[r] & 1\\
1 \ar[r] & \ar[r] K  & (B_{4}(\rp))^{(1)} \ar[r]  & \an[4] \ar[r] & 1
}}
\end{xy}
\end{equation}
The vertical arrows are inclusions, and $K'$ (resp.\ $K''$) is the kernel of the restriction of $\tau$ to $(B_{4}(\rp))^{(2)}$ (resp.\ to $(B_{4}(\rp))^{(3)}$). So $K''=(B_{4}(\rp))^{(3)}$, and since the index of $(B_{4}(\rp))^{(2)}$ (resp.\ $\Z_{2}\oplus \Z_{2}$) in $(B_{4}(\rp))^{(1)}$ (resp.\ $\an[4]$) is three, we deduce that $K'=K$, which proves part~(\ref{it:ds4a})(\ref{it:ds4aiii}). 

We now prove part~(\ref{it:ds4b}). We start by studying the quotient $(B_4(\rp))^{(1)}/(B_4(\rp))^{(3)}$. As we saw above, the elements $A_{i},B_{i},C_{i},D_{i},X_{i},Z_{i}$, $i=1,\ldots,4$ of $(B_4(\rp))^{(1)}$ are sent to the trivial element of $(B_4(\rp))^{(1)}/(B_4(\rp))^{(2)}$, and so belong to $(B_4(\rp))^{(2)}$. Hence considered as elements of $(B_4(\rp))^{(1)}/(B_4(\rp))^{(3)}$ they commute pairwise (we shall not distinguish notationally between elements of $(B_4(\rp))^{(1)}$ and their cosets in $(B_4(\rp))^{(1)}/(B_4(\rp))^{(3)}$). These were precisely the relations that we added to those of $(B_4(\rp))^{(1)}$ in order to obtain the presentation~\reqref{presd1d3} of $\Lambda$, and thus the relations of $\Lambda$ hold in $(B_4(\rp))^{(1)}/(B_4(\rp))^{(3)}$. In particular, $(B_4(\rp))^{(1)}/(B_4(\rp))^{(3)}$ is a quotient of $\Lambda$.

Since $Z_{1},Z_{3}\in (B_4(\rp))^{(2)}$, we have that $Y_{1}Z_{1}Y_{1}^{-1},Y_{1}Z_{3}Y_{1}^{-1}\in (B_4(\rp))^{(2)}$. Let $G$ denote the group obtained from $\Lambda$ by adding the following relations to the presentation~\reqref{presd1d3} of $\Lambda$:
\begin{equation}\label{eq:extra}
\left\{
\begin{aligned}
& \text{$Z_{1},Z_{3},Y_{1}Z_{1}Y_{1}^{-1},Y_{1}Z_{3}Y_{1}^{-1}$ commute pairwise}\\
& \text{and commute with $Z_{1},Z_{3},A,C$ and $D$.}
\end{aligned}\right.
\end{equation}
Once more, considered as elements of $B_{4}(\rp)$, $Z_{1},Z_{3},Y_{1}Z_{1}Y_{1}^{-1},Y_{1}Z_{3}Y_{1}^{-1},A,C$ and $D$ belong to $(B_4(\rp))^{(2)}$, and so the commutation relations of \req{extra} of $G$ also hold in $(B_4(\rp))^{(1)}/(B_4(\rp))^{(3)}$. This implies that $(B_4(\rp))^{(1)}/(B_4(\rp))^{(3)}$ is also a quotient of $G$. 

We now determine $G$ and its relationship with $(B_4(\rp))^{(1)}/(B_4(\rp))^{(3)}$. Let $L$ be the group with generators $w_{1},w_{2},w_{3},w_{4},t$ and defining relations:
\begin{equation}\label{eq:defl}
\left\{
\begin{aligned}
& \text{for all $1\leq i,j\leq 4$, $w_{i}^2=t^3=1,\, w_{i}w_{j}=w_{j}w_{i}$,} \\
& tw_{1}t^{-1}=w_{2},\, tw_{2}t^{-1}=w_{1}w_{2},\, tw_{3}t^{-1}=w_{4},\, tw_{4}t^{-1}=w_{3}w_{4}.
\end{aligned}\right.
\end{equation}
Clearly $L$ is isomorphic to $\Z_{2}^4\rtimes \Z_{3}$, where the action of conjugation by $t$ on $\ang{w_{1},\ldots,w_{4}}$ permutes cyclically the elements $w_{1},w_{2}$ and $w_{1}w_{2}$ (resp.\ $w_{3},w_{4}$ and $w_{3}w_{4}$). We define a map $\map{\psi}{L}[G]$ on the generators of $L$ as follows:
\begin{equation*}
\psi(w_{1})=Z_{1},\, \psi(w_{2})=Y_{1}Z_{1}Y_{1}^{-1},\, \psi(w_{3})=Z_{3},\, \psi(w_{4})=Y_{1}Z_{3}Y_{1}^{-1},\, \psi(t)=Y_{1}.
\end{equation*}
Since $Z_{1}^2=Z_{3}^2=1$ and $Y_{1}^3=1$ in $G$, we clearly have $(\psi(w_{i}))^2=(\psi(t))^3=1$ for $i=1,\ldots,4$. The relations~\reqref{extra} of $G$ imply that the $\psi(w_{i})$ commute pairwise. Further, $\psi(t)\psi(w_{1})(\psi(t))^{-1}=\psi(w_{2})$ by definition, and
\begin{align*}
\psi(t)\psi(w_{2})(\psi(t))^{-1}&=(\psi(t))^2\psi(w_{1})(\psi(t))^{-2}= Y_{1}^2 Z_{1} Y_{1}^{-2}=Y_{1}^{-1} Z_{1} Y_{1}^{-2}
=Z_{1}Y_{1}Z_{1}Y_{1}^{-1}\\
& =\psi(w_{1})\psi(w_{2}),
\end{align*}
using the relations $Y_{1}^3=1$ and $(Y_{1}Z_{1})^3=1$ of~\reqref{presd1d3}. Similar relations hold for $w_{3}$ and $w_{4}$, and hence $\psi$ extends to a homomorphism from $L$ to $G$. Now $\psi$ is surjective because the generating set $\brak{Z_{1},Z_{3}, Y_{1},Y_{3},A,C,D}$ of $G$ may be reduced to $\brak{Z_{1},Z_{3}, Y_{1},Y_{3}}$ using the relations~\reqref{presd1d3}. Thus $G$ is a quotient of $L$, and hence $(B_4(\rp))^{(1)}/(B_4(\rp))^{(3)}$ is also a quotient of $L$.

Let us now show that the groups $L$ and $(B_4(\rp))^{(1)}/(B_4(\rp))^{(3)}$ are isomorphic. Consider the map $\map{\phi}{(B_4(\rp))^{(1)}}[L]$ defined on the generators of $(B_4(\rp))^{(1)}$ as follows:
\begin{equation}\label{eq:defphi}
\left\{
\begin{aligned}
\phi(X_{2})&=\phi(X_{4})=\phi(A_{3})=\phi(B_{1})=\phi(B_{4})=1\\
\phi(A_{1})&=\phi(A_{4})=\phi(B_{2})=\phi(B_{3})=w_{1}w_{3}\\
\phi(C_{1})&=\phi(C_{2})=\phi(C_{3})=\phi(C_{4})=w_{1}w_{2}w_{3}w_{4}\\
\phi(D_{1})&=\phi(D_{2})=\phi(D_{3})=\phi(D_{4})=w_{2}w_{4}\\
\phi(Z_{1})&=\phi(Z_{2})=w_{1},\, \phi(Z_{3})=\phi(Z_{4})=w_{3}\\
\phi(Y_{1})&=t,\, \phi(Y_{2})=t^2,\, \phi(Y_{3})=(w_{1}w_{2}w_{3}w_{4})t^2,\, \phi(Y_{4})=(w_{1}w_{3})t.
\end{aligned}\right.
\end{equation} 
A long but straightforward calculation shows that each of the relators~\reqref{firstrel}--\reqref{lastrel} of $(B_4(\rp))^{(1)}$ is sent to the trivial element of $L$, and hence $\phi$ extends to a surjective homomorphism of $(B_4(\rp))^{(1)}$ onto $L$. Such a homomorphism sends $((B_4(\rp))^{(1)})^{(2)}=(B_4(\rp))^{(3)}$ surjectively onto $L^{(2)}$. However $L^{(2)}$ is trivial, so $\phi$ induces a surjective homomorphism $\overline{\phi}$ of $(B_4(\rp))^{(1)}/(B_4(\rp))^{(3)}$ onto $L$, and hence $L$ is a quotient of $(B_4(\rp))^{(1)}/(B_4(\rp))^{(3)}$. Since $L$ is finite and $(B_4(\rp))^{(1)}/(B_4(\rp))^{(3)}$ is a quotient of $L$ by the previous paragraph, we conclude that
\begin{equation}\label{eq:b4d1d3}
(B_4(\rp))^{(1)}/(B_4(\rp))^{(3)} \cong L\cong \Z_{2}^4 \rtimes \Z_{3},
\end{equation}
where the action is given by \req{defl}. Further, $\map{\psi}{L}[G]$ is surjective and $(B_4(\rp))^{(1)}/(B_4(\rp))^{(3)}$ is a quotient of $G$, so $G=(B_4(\rp))^{(1)}/(B_4(\rp))^{(3)}$. An easy calculation shows that $\psi^{-1}=\overline{\phi}$. From the short exact sequence
\begin{multline*}
1\to (B_4(\rp))^{(2)}/(B_4(\rp))^{(3)} \to (B_4(\rp))^{(1)}/(B_4(\rp))^{(3)} \to\\
(B_4(\rp))^{(1)}/(B_4(\rp))^{(2)}\to 1,
\end{multline*}
we see that $(B_4(\rp))^{(2)}/(B_4(\rp))^{(3)}\cong \Z_{2}^4$. It follows from the form of the isomorphism that the $\Z_{2}$-factors of $(B_4(\rp))^{(2)}/(B_4(\rp))^{(3)}$ are generated by the elements $Z_{1}$, $Z_{3}$, $Y_{1}Z_{1}Y_{1}^{-1}$ and $Y_{1}Z_{3}Y_{1}^{-1}$, and their images under $\overline{\phi}$ are $w_{1}$, $w_{3}$, $w_{2}$ and $w_{4}$ respectively.
In particular, $\overline{\phi}\left((B_4(\rp))^{(2)}/(B_4(\rp))^{(3)}\right)=\ang{w_{1},w_{2},w_{3},w_{4}}$. This completes the proof of part~(\ref{it:ds4b}).


We now prove part~(\ref{it:ds4d}). Consider the commutative diagrams~\reqref{diaggam2} and~\reqref{bigdig2}. Since $K=K'$ from part~(\ref{it:ds4a})(\ref{it:ds4aiii}) above, we have that $K\subset (B_{4}(\rp))^{(2)} \cap P_{4}(\rp)$. Conversely, since the homomorphism $(B_{4}(\rp))^{(2)}\to \Z_{2}\oplus \Z_{2}$ of \req{bigdig2} is the restriction of the permutation homomorphism $\map{\tau}{B_{4}(\rp)}[{\sn[4]}]$ to $(B_{4}(\rp))^{(2)}$, it follows that any element of $(B_{4}(\rp))^{(2)}\cap P_{4}(\rp)$ also belongs to $K$, and thus 
\begin{equation}\label{eq:kb2p4}
K= (B_{4}(\rp))^{(2)} \cap P_{4}(\rp). 
\end{equation}
Further, from the upper exact sequence of \req{bigdig2}, $(B_{4}(\rp))^{(3)}\subset K$, and since $(B_{4}(\rp))^{(3)}$ is normal in $B_{4}(\rp)$, we obtain 
\begin{equation*}
1 \to K/(B_{4}(\rp))^{(3)} \to (B_{4}(\rp))^{(2)}/(B_{4}(\rp))^{(3)} \to (B_{4}(\rp))^{(2)}/K \to 1
\end{equation*}
by taking the quotient by $(B_{4}(\rp))^{(3)}$ of the first two terms of the middle short exact sequence of \req{bigdig2}. In particular, $K/(B_{4}(\rp))^{(3)}\cong \Z_{2}\oplus \Z_{2}$, and we have a short exact sequence
\begin{equation}\label{eq:d3k}
1 \to (B_{4}(\rp))^{(3)} \to K \to \Z_{2}\oplus \Z_{2}\to 1.
\end{equation}
Recall that the $\Z_{2}$-factors of $(B_{4}(\rp))^{(2)}/(B_{4}(\rp))^{(3)}$ are generated by $Z_{1},Z_{3}, Y_{1}Z_{1}Y_{1}^{-1}$ and $Y_{1}Z_{3}Y_{1}^{-1}$. Using \req{kb2p4} and the expressions for $Z_{1},Z_{3}$ and $Y_{1}$ in terms of the standard generators of $B_{4}(\rp)$, we conclude that
\begin{equation}\label{eq:kd3quot}
K/(B_{4}(\rp))^{(3)}=\brak{1,Z_{1}Z_{3}, Y_{1}Z_{1}Z_{3}Y_{1}^{-1},Z_{1}Z_{3}Y_{1}Z_{1}Z_{3}Y_{1}^{-1}}.
\end{equation}

We now apply the Reidemeister-Schreier rewriting process to the leftmost vertical short exact sequence of~\reqref{diaggam2} to produce a set of generators of $K$. Taking $\brak{1,\rho_{1}}$ as a Schreier transversal of $\ang{\overline{\rho}}$ in $P_{4}(\rp)$ and $\set{B_{i,j},\, \rho_{k}}{1\leq i< j \leq 4,\, 1 \leq k\leq 4}$ as a generating set of $P_{4}(\rp)$, we see that the following elements constitute a generating set of $K$:
\begin{equation}\label{eq:gensk}
\left\{
\begin{aligned}
& B_{1,2}=X_{2},\, B_{1,3}=Y_{1}X_{2}Y_{1}^{-1},\, B_{1,4}= Z_{1}Y_{2}X_{2}Y_{2}^{-1}Z_{1}^{-1},\, B_{2,3}=Y_{1}Y_{2}\\
& B_{2,4}=Z_{1}Y_{2}Y_{1} Z_{1}^{-1},\, B_{3,4}=Z_{1}Z_{2},\, 
\rho_{1}B_{1,2}\rho_{1}^{-1}=A_{1}X_{4} A_{1}^{-1}\\
& \rho_{1}B_{1,3}\rho_{1}^{-1}=A_{1}Y_{4} X_{4}Y_{4}^{-1} A_{1}^{-1},\,\rho_{1}B_{1,4}\rho_{1}^{-1}=A_{1}Z_{4} Y_{3}X_{4}Y_{3}^{-1} Z_{4}^{-1} A_{1}^{-1}\\ 
& \rho_{1}^{2}=A_{1}A_{4},\,
\rho_{2}^2=B_{1}B_{4},\,
\rho_{3}^2=C_{1}C_{4},\,
\rho_{4}^2=D_{1}D_{4}\\
& \rho_{1}\rho_{2}=A_{1}B_{4},\, \rho_{1}\rho_{3}=A_{1}C_{4},\, \rho_{1}\rho_{4}=A_{1}D_{4}.
\end{aligned}\right.
\end{equation}
Note that we have also written each element in terms of the generators of the presentation of $(B_{4}(\rp))^{(1)}$ given at the beginning of the proof, we have deleted $\rho_{1}B_{2,3}\rho_{1}^{-1}$ and $\rho_{1}B_{3,4}\rho_{1}^{-1}$ from the list of generators that appear initially in the process, and that for $i=2,3,4$, we have replaced $\rho_{i}\rho_{1}^{-1}$ by $\rho_{i}^{2}=\rho_{i}\rho_{1}^{-1}\ldotp \rho_{1}\rho_{i}$.  Since $(B_{4}(\rp))^{(3)}\subset P_{4}(\rp)$, we may consider the image of $(B_{4}(\rp))^{(3)}$ in $P_{3}(\rp)$ and $P_{2}(\rp)$ under the projections 
\begin{equation*}
\text{$\map{p_{3}}{P_{4}(\rp)}[P_{3}(\rp)]$ and $\map{p_{2}}{P_{3}(\rp)}[P_{2}(\rp)]$}
\end{equation*}
obtained geometrically by forgetting the last string in each case. We claim that $p_{2}\circ p_{3}\left((B_{4}(\rp))^{(3)}\right)=\ang{\rho_{1}\rho_{2}} \cong \Z_{4}$. To see this, we first use \req{d3k} and the Reidemeister-Schreier rewriting process to obtain a generating set for $(B_{4}(\rp))^{(3)}$. This is achieved as follows. From \req{gensk}, we see that the elements of the set
\begin{equation*}
\mathcal{T}=\brak{1,\rho_{2}\rho_{1}^{-1}, \rho_{2}\rho_{1}^{-1}\ldotp \rho_{3}\rho_{1}^{-1}, \rho_{2}\rho_{1}^{-1} \ldotp\rho_{3}\rho_{1}^{-1} \ldotp\rho_{1}\rho_{2}^{-1}}
\end{equation*}
belong to $K$. Equations~\reqref{defphi} and~\reqref{kd3quot} give rise to the following commutative diagram:
\begin{equation*}
\begin{xy}*!C\xybox{%
\xymatrix{%
K/(B_{4}(\rp))^{(3)} \ar[d] \ar[r] & (B_{4}(\rp))^{(2)}/(B_{4}(\rp))^{(3)} \ar[d] \ar[r] & (B_{4}(\rp))^{(1)}/(B_{4}(\rp))^{(3)}\ar[d]\\
\ang{w_{1}w_{3}, w_{2}w_{4}} \ar[r] & \setangl{w_{i}}{i=1,\ldots,4} \ar[r] & L,}}
\end{xy}
\end{equation*}
where the horizontal arrows are inclusions, and the vertical arrows are the isomorphisms induced by $\map{\overline{\phi}}{(B_{4}(\rp))^{(1)}/(B_{4}(\rp))^{(3)}}[L]$. Equation~\reqref{gensk} yields $\rho_{2}\rho_{1}^{-1}=B_{1}A_{1}^{-1}$ and $\rho_{3}\rho_{1}^{-1}=C_{1}A_{1}^{-1}$, and considering the $K/(B_{4}(\rp))^{(3)}$-cosets of these elements and applying \req{defphi}, we obtain $\overline{\phi}(\rho_{2}\rho_{1}^{-1})=w_{1}w_{3}$ and $\overline{\phi}(\rho_{3}\rho_{1}^{-1})=w_{2}w_{4}$. 
It follows that $\mathcal{T}$ is a Schreier transversal for $K/(B_{4}(\rp))^{(3)}$ in $K$, which enables us to write down a generating set $\Sigma$ for $(B_{4}(\rp))^{(3)}$. However, to prove the claim, we do not need to study the whole list of generators. On the one hand, applying the description of the generators of $K$ given by \req{gensk}, the isomorphism $\overline{\phi}$ and \req{defphi}, we see that the elements of
\begin{equation*}
\mathcal{U}=\brak{B_{1,2}, B_{1,3}, B_{1,4}, B_{2,3}, B_{2,4}, \rho_{1}B_{1,2}\rho_{1}^{-1}, \rho_{1}B_{1,3}\rho_{1}^{-1}, \rho_{1}B_{1,4}\rho_{1}^{-1},\rho_{1}^2,\rho_{2}^2,\rho_{3}^2, \rho_{4}^2}
\end{equation*}
belong to $(B_{4}(\rp))^{(3)}$, and appear as elements of $\Sigma$. Moreover, it is clear that these elements are mapped into $\ang{\ft[2]}$ under $p_{2}\circ p_{3}$ since $\rho_{1}^2=\rho_{2}^2=\ft[2]$ in $P_{2}(\rp)$. The other elements of $\Sigma$ obtained by applying the Reidemeister-Schreier process to an element $u\in \mathcal{U}$ are just conjugates of $u$ (the conjugating elements being the non-trivial elements of $\mathcal{T}$), so also belong to $(B_{4}(\rp))^{(3)}$, and since $\ang{\ft[2]}$ is Abelian, these elements of $\Sigma$ will lie in $\ang{\ft[2]}$, which is contained in $\ang{\rho_{1}\rho_{2}}$. Hence it suffices to consider the elements of $\Sigma$ obtained by applying the Reidemeister-Schreier rewriting process to the three remaining elements $\rho_{1}\rho_{i}$, $i=2,3,4$, of \req{gensk}. To do this, note that under  identification of $(B_{4}(\rp))^{(1)}/(B_{4}(\rp))^{(3)}$ with $L$, the elements of $\mathcal{T}$ project respectively to $1$, $w_{1}w_{3}$, $w_{1}w_{2}w_{3}w_{4}$ and $w_{2}w_{4}$, while $\rho_{1}\rho_{2}$ projects to $w_{1}w_{3}$, $\rho_{1}\rho_{3}$ projects to $w_{2}w_{4}$, and $\rho_{1}\rho_{4}$ projects to $w_{1}w_{2}w_{3}w_{4}$. The non-trivial elements of $\Sigma$ arising as conjugates of $\rho_{1}\rho_{i}$ are as follows:
\begin{enumerate}
\item $i=2$: $\rho_{1}\rho_{2}\rho_{1}\rho_{2}^{-1}$, $\rho_{2}^2$, $\rho_{2}\rho_{1}^{-1} \rho_{3}\rho_{2}^2 \rho_{3}^{-1} \rho_{1}\rho_{2}^{-1}$, $\rho_{2}\rho_{1}^{-1} \rho_{3}\rho_{2}^{-1} \rho_{1}\rho_{2} \rho_{1}\rho_{3}^{-1} \rho_{1}\rho_{2}^{-1}$.
\item $i=3$: $\rho_{1} \rho_{3} \rho_{2}\rho_{3}^{-1} \rho_{1}\rho_{2}^{-1}$, $\rho_{2} \rho_{3} \rho_{1}\rho_{3}^{-1} \rho_{1}\rho_{2}^{-1}$, $\rho_{2}\rho_{1}^{-1} \rho_{3}^2 \rho_{1}\rho_{2}^{-1}$, $\rho_{2}\rho_{1}^{-1} \rho_{3}\rho_{2}^{-1}\rho_{1}\rho_{3}$.
\item $i=4$: $\rho_{1}\rho_{4} \rho_{1}\rho_{3}^{-1} \rho_{1}\rho_{2}^{-1}$, $\rho_{2}\rho_{4} \rho_{2}\rho_{3}^{-1} \rho_{1}\rho_{2}^{-1}$, $\rho_{2}\rho_{1}^{-1} \rho_{3}\rho_{4}$, $\rho_{2}\rho_{1}^{-1} \rho_{3}\rho_{2}^{-1} \rho_{1}\rho_{4} \rho_{1}\rho_{2}^{-1}$.
\end{enumerate}
Under the projection $p_{2}\circ p_{3}$, the elements for the cases $i=2,3$ project to elements of $\ang{\ft[2]}$, while those for the case $i=4$ project to $\rho_{1}\rho_{2}$ or its inverse. We conclude that $p_{2}\circ p_{3}\left((B_{4}(\rp))^{(3)}\right)=\ang{\rho_{1}\rho_{2}} \cong \Z_{4}$, which proves the claim. Thus the restriction
\begin{equation*}
\map{p_{2}\left\lvert_{p_{3}\left((B_{4}(\rp))^{(3)}\right)}\right.}{p_{3}((B_{4}(\rp))^{(3)})}[\ang{\rho_{1}\rho_{2}}]
\end{equation*}
of $p_{2}$ to $p_{3} \left((B_{4}(\rp))^{(3)} \right)$ is surjective. 

Now consider the following commutative diagram of short exact sequences:
\begin{equation}\label{eq:indfour}
\begin{xy}*!C\xybox{%
\xymatrix{%
1 \ar[r] &  \ker{p_{2}\left\lvert_{p_{3} \left((B_{4}(\rp))^{(3)} \right)}\right.} \ar[r] \ar[d] & p_{3} \left((B_{4}(\rp))^{(3)} \right) \ar[d] \ar[rrr]^>>>>>>>>>>>>>>>>>{p_{2}\left\lvert_{p_{3} \left((B_{4}(\rp))^{(3)} \right)}\right.} & & & \ang{\rho_{1}\rho_{2}} \ar[d]\ar[r] & 1\\
1 \ar[r] & \ar[r] \F[2](\rho_{3},B_{2,3})  & P_{3}(\rp)  \ar[rrr]^{p_{2}} & & & P_{2}(\rp) \ar[r] & 1.}}
\end{xy}
\end{equation}
The lower short exact sequence is that of \req{split} with $m=2$ and $n=1$ (here $p_{\ast}=p_{2}$), while the vertical arrows are inclusions. It follows that
\begin{equation}\label{eq:kerp2inter}
\ker{p_{2}\left\lvert_{p_{3} \left((B_{4}(\rp))^{(3)} \right)}\right.}= p_{3} \left((B_{4}(\rp))^{(3)} \right)\cap \F[2](\rho_{3},B_{2,3}).
\end{equation}
Since $K\subset P_{4}(\rp)$, $p_{3}$ restricts to $K$, and we have the following commutative diagram:
\begin{equation}\label{eq:commk}
\begin{xy}*!C\xybox{%
\xymatrix{%
(B_{4}(\rp))^{(3)} \ar[rrr]^>>>>>>>>>>>>>>>{p_{3}\left\lvert_{(B_{4}(\rp))^{(3)}}\right.} \ar[d] & & & p_{3}\left((B_{4}(\rp))^{(3)}\right) \ar[d]\\
K \ar[rrr]^{p_{3}\left\lvert_{K}\right.} & & & P_{3}(\rp).
}}
\end{xy}
\end{equation}
Again the vertical arrows are inclusions. Considered as elements of $P_{4}(\rp)$, $B_{i,j}$, $1\leq i<j\leq 3$ and $\rho_{k}\rho_{4}$, $1\leq k\leq 3$, belong to $K$ by \req{diaggam2}, and we deduce that the restriction of $p_{3}$ to $K$ is surjective. Since $(B_{4}(\rp))^{(3)}$ is of index four in $K$, we conclude that $p_{3}\left((B_{4}(\rp))^{(3)}\right)$ is of index at most four in $P_{3}(\rp)$. 

Conversely, consider the Abelianisation of $P_{3}(\rp)$. From \req{vanbuskirk} and the action of $\quat$ on $\F[2](\rho_{3},B_{2,3})$ described by equations~\reqref{actr12}, \reqref{conjrho3} and~\reqref{conjrho3a}, we see that $\bnab{3}\cong \Z_{2}^3$, and that the Abelianisation homomorphism $\map{\pi}{P_{3}(\rp)}[\Z_{2}^3]$ sends each of $\rho_{i}$, $i=1,2,3$, to a distinct $\Z_{2}$-factor, and the $B_{i,j}$, $1\leq i<j\leq 3$, to the trivial element. Under $p_{3}$, the elements of $\Sigma$ are sent to the trivial element of $\Z_{2}^3$, with the exception of those elements obtained via the Reidemeister-Schreier rewriting process using $\rho_{1}\rho_{4}$, which are sent to $(\overline{1}, \overline{1}, \overline{1})$. It follows that
\begin{equation*}
\pi \left(p_{3}\left((B_{4}(\rp))^{(3)}\right)\right)=\ang{(\overline{1}, \overline{1}, \overline{1})}\cong \Z_{2},
\end{equation*}
and so $\pi\left(p_{3}((B_{4}(\rp))^{(3)})\right)$ is of index four in $\Z_{2}^3$. We conclude from the following commutative diagram: 
\begin{equation*}
\xymatrix{%
p_{3}\left((B_{4}(\rp))^{(3)}\right) \ar[r] \ar[d]_{\pi} & P_{3}(\rp) \ar[d]_{\pi}\\
\pi\left(p_{3}\left((B_{4}(\rp))^{(3)}\right)\right) \ar[r]  & \Z_{2}^3,
}
\end{equation*}
whose horizontal arrows are inclusions and whose vertical arrows are surjections, that $p_{3}\left((B_{4}(\rp))^{(3)}\right)$ is of index at least four in $P_{3}(\rp)$. From the previous paragraph, we conclude this index is exactly four, and since $\ang{\rho_{1}\rho_{2}}$ is of index two in $P_{2}(\rp)$, it follows from equations~\reqref{indfour} and~\reqref{kerp2inter} that $p_{3}\left((B_{4}(\rp))^{(3)}\right)\cap \F[2](\rho_{3},B_{2,3})$ is of index two in $\F[2](\rho_{3},B_{2,3})$. Since $B_{2,3}\in \Sigma$ (as an element of $P_{4}(\rp)$), we have that $B_{2,3}\in p_{3}\left((B_{4}(\rp))^{(3)}\right)$ (as an element of $P_{3}(\rp)$). Thus under the canonical homomorphism 
\begin{equation*}
\F[2](\rho_{3},B_{2,3})\to \F[2](\rho_{3},B_{2,3})\left/\left(p_{3}\left((B_{4}(\rp))^{(3)}\right)\cap \F[2](\rho_{3},B_{2,3})\right)\right.\cong\Z_{2},
\end{equation*}
$B_{2,3}$ is sent to $\overline{0}$, so $\rho_{3}$ must be sent to $\overline{1}$, and hence the kernel of this homomorphism is given by:
\begin{equation}\label{eq:kerprojf2}
p_{3}\left((B_{4}(\rp))^{(3)}\right)\cap \F[2](\rho_{3},B_{2,3})=\F[3](B_{2,3},\rho_{3}^2, \rho_{3}B_{2,3}\rho_{3}^{-1}).
\end{equation}
From \req{gensk}, $\rho_{4}\rho_{3}\rho_{2}\rho_{1}=D_{1}C_{4} B_{1} A_{4}$, and using \req{defphi} and the isomorphism of \req{b4d1d3}, we see that $\rho_{4}\rho_{3}\rho_{2}\rho_{1}\in (B_{4}(\rp))^{(3)}$, and so $\rho_{3}\rho_{2}\rho_{1}\in p_{3}\left((B_{4}(\rp))^{(3)}\right)$ and $\rho_{2}\rho_{1}\in p_{2}\circ p_{3}\left((B_{4}(\rp))^{(3)}\right)$. But we know that each of these three elements is of order four in its respective group~\cite[Proposition~26 and Remark~27]{GG4}, and hence it follows from the upper sequence of \req{indfour} and \req{kerprojf2} that 
\begin{equation}\label{eq:f3z4}
p_{3}\left((B_{4}(\rp))^{(3)}\right)\cong \F[3] \rtimes \Z_{4}.
\end{equation}
We shall determine the action shortly. Returning to \req{commk}, both of the horizontal restrictions $p_{3}\left\lvert_{(B_{4}(\rp))^{(3)}}\right.$  and $p_{3}\left\lvert_{K}\right.$ are surjective, and since $(B_{4}(\rp))^{(3)}$ (resp.\ $p_{3}((B_{4}(\rp))^{(3)})$) is of index four in $K$ (resp.\ $P_{3}(\rp)$), we obtain $\ker{p_{3}\left\lvert_{(B_{4}(\rp))^{(3)}}\right.}=\ker{p_{3}\left\lvert_{K}\right.}$. Thus from the upper homomorphism of \req{commk} and \req{f3z4}, we have a short exact sequence
\begin{equation*}
1 \to \ker{p_{3}\left\lvert_{K}\right.} \to (B_{4}(\rp))^{(3)} 
\xrightarrow{p_{3}\left\lvert_{(B_{4}(\rp))^{(3)}}\right.} 
\F[3] \rtimes \Z_{4}\to 1.
\end{equation*}
From equations~\reqref{indfour}, \reqref{kerp2inter} and~\reqref{kerprojf2}, a basis of the $\F[3]$-factor of the quotient is given by $\brak{B_{2,3},\rho_{3}^2,\rho_{3} B_{2,3}\rho_{3}^{-1}}$, and by \reqref{indfour} and the above discussion, we may take $a^3=\rho_{3}\rho_{2}\rho_{1}$ to be a generator of the $\Z_{4}$-factor. Using \req{cyclicperm2}, we see that the action of $\Z_{4}$ on $\F[3]$ is given by 
\begin{equation}\label{eq:conja3}
\left\{
\begin{aligned}
a^3 B_{2,3}a^{-3}&= B_{2,3}^{-1}\\
a^3\rho_{3}^2 a^{-3}&= \rho_{3}^{-2}\\
a^3\rho_{3} B_{2,3}\rho_{3}^{-1}a^{-3}&=\rho_{3}^{-2} \ldotp \rho_{3} B_{2,3}^{-1} \rho_{3}^{-1}\ldotp \rho_{3}^{2}.
\end{aligned}
\right.
\end{equation}
Consider the map $\map{s}{\F[3] \rtimes \Z_{4}}[(B_{4}(\rp))^{(3)}]$ defined on the generators of $\F[3] \rtimes \Z_{4}$ by:
\begin{equation*}
\left\{
\begin{aligned}
x &\mapsto x \quad \text{for $x\in \brak{B_{2,3},\rho_{3}^2,\rho_{3} B_{2,3}\rho_{3}^{-1}}$}\\
a^3 & \mapsto a^4.
\end{aligned}
\right.
\end{equation*}
Note that the elements on the right hand-side are considered to be elements in $B_{4}(\rp)$. Using the given generating set of $K$, we have
\begin{equation*}
B_{2,3}=Y_{1}Y_{2},\; \rho_{3}^2=C_{1}C_{4},\; \rho_{3} B_{2,3}\rho_{3}^{-1}= C_{1}Y_{4}Y_{3}C_{1}^{-1},\; a^4=\rho_{4}\rho_{3}\rho_{2}\rho_{1}= D_{1}C_{4}B_{1}A_{4}.
\end{equation*}
By equations~\reqref{defphi} and~\reqref{b4d1d3}, we see that these elements belong to $(B_{4}(\rp))^{(3)}$, so the map $s$ is well defined. Using \req{cyclicperm2} once more, we see that the action of $a^{4}=s(a^{3})$ (which is of order $4$) on $s(x)$, $x\in \brak{B_{2,3},\rho_{3}^2,\rho_{3} B_{2,3}\rho_{3}^{-1}}$, is also given by \req{conja3}, up to replacing $a^{3}$ by $a^{4}$. This shows that $s$ extends to a homomorphism from $\F[3] \rtimes \Z_{4}$ to $(B_{4}(\rp))^{(3)}$. It is then clear that $s$ is a section for $p_{3}\left\lvert_{(B_{4}(\rp))^{(3)}}\right.$, and hence
\begin{equation*}
(B_{4}(\rp))^{(3)}\cong \ker{p_{3}\left\lvert_{K}\right.} \rtimes (\F[3] \rtimes \Z_{4}).
\end{equation*}

From the following commutative diagram of short exact sequences,
\begin{equation*}
\xymatrix{%
& 1 \ar[d] & 1 \ar[d] & & \\
1 \ar[r] & \ker{p_{3}\left\lvert_{K}\right.}\cong\F[5] \ar[r] \ar[d] & \ker{p_{3}}\cong\F[3] \ar[r]^>>>>>{\alpha'\left\lvert_{\F[3]}\right.} \ar[d] & \ang{\overline{\rho}} \ar[r] \ar@{=}[d] & 1\\
1 \ar[r] & K \ar[r] \ar[d]^{p_{3}\left\lvert_{K}\right.} & P_{4}(\rp) \ar[r]^>>>>>>>{\alpha'} \ar[d]^{p_{3}} & \ang{\overline{\rho}} \ar[r]  & 1\\
& P_{3}(\rp) \ar[d] \ar@{=}[r] & P_{3}(\rp) \ar[d] & &\\
& 1 & 1 & &
}
\end{equation*}
we see that $\ker{p_{3}\left\lvert_{K}\right.}$ is also the kernel of the restriction of $\alpha'$ to $\ker{p_{3}}$ which is the free subgroup of $P_{4}(\rp)$ of rank three with basis $\brak{B_{1,4}, B_{2,4}, \rho_{4}}$. It thus follows that $\ker{p_{3}\left\lvert_{K}\right.}$ is a free group of rank five with basis $\brak{B_{1,4}, B_{2,4}, \rho_{4}^2,\rho_{4}B_{1,4}\rho_{4}^{-1}, \rho_{4}B_{2,4}\rho_{4}^{-1}}$. We conclude that 
\begin{align}
(B_{4}(\rp))^{(3)}& \cong \F[5](B_{1,4}, B_{2,4}, \rho_{4}^2,\rho_{4}B_{1,4}\rho_{4}^{-1}, \rho_{4}B_{2,4}\rho_{4}^{-1})\rtimes \left(\F[3](B_{2,3},\rho_{3}^2, \rho_{3}B_{2,3}\rho_{3}^{-1}) \rtimes \Z_{4}\right)\notag\\
& \cong \left(\F[5](B_{1,4}, B_{2,4}, \rho_{4}^2,\rho_{4}B_{1,4}\rho_{4}^{-1}, \rho_{4}B_{2,4}\rho_{4}^{-1})\rtimes \F[3](B_{2,3},\rho_{3}^2, \rho_{3}B_{2,3}\rho_{3}^{-1})\right) \rtimes \Z_{4}.\label{eq:d3b4}
\end{align}
As we already mentioned above, $\rho_{4}\rho_{3}\rho_{2}\rho_{1}$ belongs to $(B_{4}(\rp))^{(3)}$, and it projects to the generator of the $\Z_{4}$-factor of $p_{3}((B_{4}(\rp))^{(3)})$, so may be taken as a generator of the $\Z_{4}$-factor in \req{d3b4}. To determine completely $(B_{4}(\rp))^{(3)}$, it just remains to calculate the actions. By \req{cyclicperm2}, the action of $\Z_{4}$ on the given generators of $\F[5]\rtimes \F[3]$ is:
\begin{equation}\label{eq:actz4}
\left\{
\begin{aligned}
B_{1,4} & \mapsto \rho_{4}^2\ldotp  B_{1,4}^{-1} \ldotp\rho_{4}^{-2} &
B_{2,4}& \mapsto \rho_{4}^2\ldotp  B_{1,4} \ldotp B_{2,4}^{-1} \ldotp B_{1,4}^{-1} \ldotp\rho_{4}^{-2}\\
\rho_{4}^2 & \mapsto \rho_{4}^{-2} & 
\rho_{4}B_{1,4}\rho_{4}^{-1}& \mapsto \rho_{4}B_{1,4}^{-1}\rho_{4}^{-1}\\
\rho_{4}B_{2,4}\rho_{4}^{-1}& \mapsto  \rho_{4} B_{1,4} \rho_{4}^{-1}\ldotp \rho_{4} B_{2,4}^{-1} \rho_{4}^{-1}\ldotp \rho_{4}B_{1,4}^{-1}\rho_{4}^{-1} &
B_{2,3}& \mapsto B_{2,3}^{-1}\\
\rho_{3}^2 & \mapsto \rho_{3}^{-2} &
\rho_{3}B_{2,3}\rho_{3}^{-1}&\mapsto \rho_{3}^{-2}\ldotp \rho_{3}B_{2,3}^{-1}\rho_{3}^{-1}\ldotp \rho_{3}^{2},
\end{aligned}\right.
\end{equation}
and so the action of $\Z_{4}$ by conjugation on the Abelianisation of $\F[5]\rtimes \F[3]$ is $-\id$. As for the action by conjugation of $\F[3]$ on $\F[5]$, we have:
\begin{equation}\label{eq:actf3f5a}
B_{2,3}:
 \left\{
\begin{aligned}
B_{1,4}& \mapsto B_{1,4}\\
B_{2,4}& \mapsto \rho_{4}^2 \ldotp B_{1,4} \ldotp B_{2,4} \ldotp B_{1,4}^{-1} \ldotp\rho_{4}^{-2}\\
\rho_{4}^2 & \mapsto \rho_{4}^2\\
\rho_{4}B_{1,4}\rho_{4}^{-1}& \mapsto \rho_{4}B_{1,4}\rho_{4}^{-1}\\
\rho_{4}B_{2,4}\rho_{4}^{-1}& \mapsto \rho_{4}^2 \cdot \rho_{4}B_{1,4}\rho_{4}^{-1}\cdot \rho_{4}B_{2,4}\rho_{4}^{-1} \cdot \rho_{4}B_{1,4}^{-1}\rho_{4}^{-1} \cdot \rho_{4}^{-2}
\end{aligned}\right.
\end{equation}
\begin{equation}\label{eq:actf3f5b}
\rho_{3}B_{2,3}\rho_{3}^{-1}:
 \left\{
\begin{aligned}
B_{1,4} \mapsto &
B_{2,4}^{-1} \cdot \rho_{4}^{-2} \cdot \rho_{4}B_{2,4}^{-1}\rho_{4}^{-1} \cdot \rho_{4} B_{1,4}^{-1}\rho_{4}^{-1} \cdot  
B_{2,4} \cdot \rho_{4}B_{1,4} \rho_{4}^{-1} \cdot \\
& \rho_{4} B_{2,4} \rho_{4}^{-1} \cdot \rho_{4}^2 \cdot B_{1,4} \cdot \rho_{4}^{-2} \cdot \rho_{4}B_{2,4}^{-1}\rho_{4}^{-1}\cdot \rho_{4}B_{1,4}^{-1}\rho_{4}^{-1}\cdot B_{2,4}^{-1} \cdot\\
& \rho_{4}B_{1,4} \rho_{4}^{-1} \cdot \rho_{4}B_{2,4}\rho_{4}^{-1}\cdot \rho_{4}^2 \cdot B_{2,4}\\
B_{2,4}\mapsto & 
B_{2,4}^{-1} \cdot \rho_{4}^{-2} \cdot \rho_{4}B_{2,4}^{-1}\rho_{4}^{-1} \cdot \rho_{4}B_{1,4}^{-1}\rho_{4}^{-1} \cdot B_{2,4} \cdot \rho_{4}B_{1,4}\rho_{4}^{-1} \cdot \\
& \rho_{4}B_{2,4}\rho_{4}^{-1}\cdot\rho_{4}^2 \cdot B_{2,4}\\
\rho_{4}^2 \mapsto & 
B_{2,4}^{-1} \cdot B_{1,4}^{-1} \cdot \rho_{4}^{-2}\cdot \rho_{4}B_{2,4}^{-1}\rho_{4}^{-1} \cdot \rho_{4}^2 \cdot B_{1,4} \cdot B_{2,4} \cdot \rho_{4}B_{1,4}^{-1}\rho_{4}^{-1} \cdot \\
& B_{2,4}^{-1} \cdot \rho_{4}B_{1,4}\rho_{4}^{-1} \cdot \rho_{4}B_{2,4}\rho_{4}^{-1}\cdot \rho_{4}^2 \cdot B_{2,4}\\
\rho_{4}B_{1,4}\rho_{4}^{-1}\mapsto & \rho_{4}B_{1,4}\rho_{4}^{-1}\\
\rho_{4}B_{2,4}\rho_{4}^{-1}\mapsto &
B_{2,4}^{-1} \cdot B_{1,4}^{-1} \cdot \rho_{4}^{-2}\cdot \rho_{4}B_{2,4}\rho_{4}^{-1} \cdot 
\rho_{4}^2 \cdot B_{1,4} \cdot B_{2,4}\end{aligned}\right.
\end{equation}
\begin{equation}\label{eq:actf3f5c}
\rho_{3}^2:
\left\{
\begin{aligned}
B_{1,4} \mapsto &\rho_{4}B_{1,4}\rho_{4}^{-1} \cdot \rho_{4}B_{2,4}\rho_{4}^{-1} \cdot B_{2,4}^{-1} \cdot B_{1,4} \cdot B_{2,4} \cdot \rho_{4}B_{2,4}^{-1}\rho_{4}^{-1} \cdot\\
& \rho_{4}B_{1,4}^{-1}\rho_{4}^{-1}\\
B_{2,4} \mapsto &\rho_{4}B_{1,4}\rho_{4}^{-1} \cdot \rho_{4}B_{2,4}\rho_{4}^{-1} \cdot B_{2,4}^{-1} \cdot B_{1,4}^{-1} \cdot B_{2,4}\cdot B_{1,4} \cdot B_{2,4} \cdot \\
& \rho_{4}B_{2,4}^{-1}\rho_{4}^{-1} \cdot \rho_{4}B_{1,4}^{-1}\rho_{4}^{-1}\\
\rho_{4}^2  \mapsto &\rho_{4}B_{1,4}\rho_{4}^{-1} \cdot \rho_{4}B_{2,4}\rho_{4}^{-1} \cdot \rho_{4}^2 \cdot \rho_{4}B_{2,4}^{-1}\rho_{4}^{-1} \cdot \rho_{4}B_{1,4}^{-1}\rho_{4}^{-1}\\
\rho_{4}B_{1,4}\rho_{4}^{-1} \mapsto &\rho_{4}B_{1,4}\rho_{4}^{-1}\\
\rho_{4}B_{2,4}\rho_{4}^{-1} \mapsto &\rho_{4}B_{2,4}\rho_{4}^{-1},
\end{aligned}\right.
\end{equation}
using the relation $B_{1,4}B_{2,4}B_{3,4}=\rho_{4}^{-2}$ in $P_{4}(\rp)$. 
In all cases, the action of the given generators of $\F[3]$ on the Abelianisation of $\F[5]$ is trivial. We thus conclude that
\begin{equation*}
\gpab[\left((B_{4}(\rp))^{(3)}\right)]= (B_{4}(\rp))^{(3)}/(B_{4}(\rp))^{(4)}\cong \Z_{2}^8 \oplus \Z_{4},
\end{equation*}
the $\Z_{2}$-factors arising from the fact that the action of $\Z_{4}$ on $\F[5]\rtimes \F[3]$ is $-\id$. Consider the following short exact sequence:
\begin{equation*}
1 \to (B_{4}(\rp))^{(4)} \to (B_{4}(\rp))^{(3)} \to (B_{4}(\rp))^{(3)}/(B_{4}(\rp))^{(4)} \to 1.
\end{equation*}
The $\Z_{4}$-factor of $(B_{4}(\rp))^{(3)}$ is mapped bijectively onto the $\Z_{4}$-factor of the quotient $(B_{4}(\rp))^{(3)}/(B_{4}(\rp))^{(4)}$, so the kernel $(B_{4}(\rp))^{(4)}$ of the projection $(B_{4}(\rp))^{(3)} \to (B_{4}(\rp))^{(3)}/(B_{4}(\rp))^{(4)}$ is the restriction of this projection to $\F[5]\rtimes \F[3]$. From the form of the action of $\F[3]$ on $\F[5]$, this restriction is the composition of the Abelianisation  $\F[5]\rtimes \F[3] \to \Z^5 \oplus \Z^3$, followed by the homomorphism $\Z^5 \oplus \Z^3 \to \Z_{2}^5 \oplus \Z_{2}^3$ which takes the coordinates modulo $2$. We see that
\begin{equation}\label{eq:act129}
(B_{4}(\rp))^{(4)}\cong \F[129]\rtimes \F[17],
\end{equation}
where $\F[129]$ (resp.\ $\F[17]$) is the kernel of the restriction $\F[5]\to \Z_{2}^5$ (resp.\ of $\F[3]\to \Z_{2}^3$) of this composition to the first (resp.\ second) factor, and the action is that induced by that of $\F[3]$ on $\F[5]$. It is then clear that for all $i\geq 0$, $(B_{4}(\rp))^{(4+i)}\cong (\F[129]\rtimes \F[17])^{(i)}$. This completes the proof of  part~(\ref{it:ds4d}), and thus that of \reth{dsb4}.
\end{proof}

\begin{rem}
In order to decide whether $B_{4}(\rp)$ is residually soluble, it would be useful to know the form of the action in \req{act129}. If the product $\F[129]\rtimes \F[17]$ were almost direct (\emph{i.e.}\ the action of $\F[17]$ on the Abelianisation of $\F[129]$ were trivial) then $B_{4}(\rp)$ would be residually soluble~\cite{FR1}. However, this is not the case. To see this, first consider the following basis $(e_{1}, \ldots, e_{5})$ of $\F[5]$:
\begin{equation}\label{eq:eibasis}
\begin{aligned}
&\text{$e_{1}=B_{1,4}$, $e_{2}=B_{2,4}$, $e_{3}=\rho_{4}B_{1,4}\rho_{4}^{-1} \ldotp \rho_{4}B_{2,4}\rho_{4}^{-1} \ldotp \rho_{4}^2$,}\\
&\text{$e_{4}=\rho_{4}B_{1,4}\rho_{4}^{-1}$ and $e_{5}=\rho_{4}B_{2,4}\rho_{4}^{-1}$.}
\end{aligned}
\end{equation}
From \req{actf3f5c}, the action of $\rho_{3}^2$ by conjugation on this basis is given by the following automorphism of $\F[5]$:
\begin{equation}\label{eq:actphi3}
\phi_{\rho_{3}^2}: 
\left\{
\begin{aligned}
e_{1} & \mapsto e_{4}e_{5} e_{2}^{-1} e_{1} e_{2} e_{5}^{-1} e_{4}^{-1}\\
e_{2} & \mapsto e_{4}e_{5} e_{2}^{-1} e_{1}^{-1} e_{2} e_{1} e_{2} e_{5}^{-1} e_{4}^{-1}\\
e_{3} & \mapsto e_{4}e_{5} e_{3} e_{5}^{-1} e_{4}^{-1}\\
e_{4} & \mapsto e_{4}\\
e_{5} & \mapsto e_{5}.
\end{aligned}
\right.
\end{equation}
It follows from the form of the projection $\F[3]\to \Z_{2}^3$ that $\rho_{3}^4$ belongs to the kernel $\F[17]$. We will calculate the action of the corresponding automorphism $\phi_{\rho_{3}^4}$ on a certain element of $\F[129]$. To do this, we first determine a basis of $\F[129]$ using the Reidemeister-Schreier rewriting process. A suitable transversal for the kernel of $\F[5]\to \Z_{2}^5$ relative to the basis of \req{eibasis} is the word
\begin{equation*}
\tau=e_{1}e_{2}e_{1} e_{3} e_{1}e_{2}e_{1} e_{4} e_{1}e_{2}e_{1} e_{3} e_{1}e_{2}e_{1} e_{5} e_{1}e_{2}e_{1} e_{3} e_{1}e_{2}e_{1} e_{4} e_{1}e_{2}e_{1} e_{3} e_{1}e_{2}e_{1}.
\end{equation*}
Let $\tau_{0}$ denote the empty word, for $i=1,\ldots, 31$, let $\tau_{i}$ be the subword of $\tau$ consisting of the first $i$ letters, and let $\overline{w}$ denote the Schreier representative of the word $w=w(e_{1},\ldots,e_{5})$. Deleting the trivial elements that appear in the set
\begin{equation*}
\set{\tau_{i}e_{j}
(\overline{\tau_{i}e_{j}})^{-1}}{0\leq i\leq 31,\, 1\leq j\leq 5},
\end{equation*}
gives rise to a basis of $\F[129]$, and thus a basis for the Abelianisation $\Z_{129}$ of $\F[129]$ (we shall not distinguish notationally between a basis element of $\F[129]$ and its projection in $\Z_{129}$). Using \req{actphi3}, a long but straightforward calculation in $\Z_{129}$ shows that
\begin{align*}
\phi_{\rho_{3}^4}(\tau_{5}e_{3}\tau_{2}^{-1})=& e_{2}^{-1} \tau_{3}^{-1} \ldotp \tau_{3}e_{1} \tau_{2}^{-1} \ldotp \tau_{2}e_{2}\tau_{1}^{-1} \ldotp \tau_{1}e_{1} \ldotp e_{2}\tau_{3}^{-1} \ldotp \tau_{4}e_{2}^{-1}\tau_{7}^{-1} \ldotp\tau_{5}e_{1}\tau_{4}^{-1} \ldotp \tau_{4}e_{2}\tau_{7}^{-1} \ldotp\\
&  \tau_{7}e_{1}\tau_{6}^{-1} \ldotp \tau_{6}e_{2}\tau_{5}^{-1} \ldotp \tau_{5}e_{3}\tau_{2}^{-1}.
\end{align*}
Each of the terms appearing on the right hand-side of this equality, as well as $\tau_{5}e_{3}\tau_{2}^{-1}$, belongs to the given basis of $\Z_{129}$, and so the induced action of $\F[17]$ on $\Z_{129}$ is non trivial. This proves that the semi-direct product $\F[129]\rtimes \F[17]$ is not almost direct. It thus remains an open question as to whether $B_{4}(\rp)$ is residually soluble.
\end{rem}

\section{A presentation of $\Gamma_2(B_n(\rp))$, $n\geq 3$}\label{sec:prescom}

In this section, we derive a presentation of $\Gamma_2(B_n(\rp))$ obtained using the Reidemeister-Schreier rewriting process. 

\begin{prop}\label{prop:fullpres}
Let $n\geq 3$. The following constitutes a presentation of the group $\Gamma_2(B_n(\rp))$:
\begin{enumerate}
\item[\underline{\textbf{generators:}}] 
\begin{align*}
\alpha_{i}&= \sigma_{i}\sigma_{1}^{-1},  
\quad \gamma_{i}= \sigma_{1}\rho_{1} \sigma_{i}\sigma_{1}^{-1} \rho_{1}^{-1}\sigma_{1}^{-1}  &&\text{for $i=2,\ldots,n-1$}\\
\beta_{i}&= \sigma_{1}\sigma_{i}, 
\quad \tau_{i}= \sigma_{1}\rho_{1} \sigma_{1}\sigma_{i} \rho_{1}^{-1}\sigma_{1}^{-1} &&\text{for $i=1,\ldots,n-1$}\\
\eta_{j}&= \rho_{j}\sigma_{1}^{-1}\rho_{1}^{-1}\sigma_{1}^{-1},  \quad \theta_{j}= \sigma_{1}\rho_{1} \rho_{j}\sigma_{1}^{-1},  \quad \lambda_{j}= \sigma_{1}\rho_{1} \sigma_{1}\rho_{j} &&\text{for $j=1,\ldots,n$}\\
\kappa_{j}&= \sigma_{1}\rho_{j}\rho_{1}^{-1}\sigma_{1}^{-1} &&\text{for $j=2,\ldots,n$.}
\end{align*}
To simplify the expression of the relators, we set $\alpha_{1}=\gamma_{1}=\kappa_{1}=1$.
\item[\underline{\textbf{relators:}}] \mbox{}
\begin{enumerate}[(a)]
\item For all $1\leq i,j\leq n-1$, $\lvert i-j\rvert \geq 2$, 
\begin{equation*}
\alpha_{i} \beta_{j} \beta_{i}^{-1}\alpha_{j}^{-1}, \quad \beta_{i} \alpha_{j} \alpha_{i}^{-1}\beta_{j}^{-1}, \quad \gamma_{i} \tau_{j} \tau_{i}^{-1}\gamma_{j}^{-1}, \quad \tau_{i} \gamma_{j} \gamma_{i}^{-1}\tau_{j}^{-1}.
\end{equation*}
\item For all $1\leq i\leq n-2$,
\begin{align*}
& \alpha_{i}\beta_{i+1}\alpha_{i}\alpha_{i+1}^{-1}\beta_{i}^{-1}\alpha_{i+1}^{-1}, 
\qquad \beta_{i} \alpha_{i+1} \beta_{i} \beta_{i+1}^{-1} \alpha_{i}^{-1} \beta_{i+1}^{-1},\\
& \gamma_{i} \tau_{i+1} \gamma_{i} \gamma_{i+1}^{-1} \tau_{i}^{-1} \gamma_{i+1}^{-1}, \qquad \tau_{i} \gamma_{i+1} \tau_{i} \tau_{i+1}^{-1} \gamma_{i}^{-1} \tau_{i+1}^{-1}.
\end{align*}
\item For all $1\leq i\leq n-1$ and $1\leq j\leq n$ with $j\neq i,i+1$,
\begin{equation*}
\alpha_{i}\kappa_{j}\tau_{i}^{-1}\eta_{j}^{-1}, \quad \beta_{i}\eta_{j}\gamma_{i}^{-1}\kappa_{j}^{-1}, \quad \gamma_{i}\lambda_{j}\beta_{i}^{-1}\theta_{j}^{-1}, \quad 
\tau_{i}\theta_{j}\alpha_{i}^{-1}\lambda_{j}^{-1}.
\end{equation*}

\item For all $1\leq i\leq n-1$,
\begin{equation*}
\beta_{i}^{-1} \kappa_{i} \tau_{i}^{-1} \eta_{i+1}^{-1}, \quad \alpha_{i}^{-1} \eta_{i} \gamma_{i}^{-1} \kappa_{i+1}^{-1}, \quad \tau_{i}^{-1} \lambda_{i} \beta_{i}^{-1} \theta_{i+1}^{-1}, \quad \gamma_{i}^{-1} \theta_{i} \alpha_{i}^{-1} \lambda_{i+1}^{-1}.
\end{equation*}

\item For all $1\leq i\leq n-1$,
\begin{align*}
&\lambda_{i+1}^{-1} \eta_{i}^{-1} \eta_{i+1} \lambda_{i}\beta_{i}^{-1} \alpha_{i}^{-1},\quad 
\theta_{i+1}^{-1} \kappa_{i}^{-1} \kappa_{i+1} \theta_{i} \alpha_{i}^{-1}\beta_{i}^{-1},\\
&\kappa_{i+1}^{-1} \theta_{i}^{-1} \theta_{i+1} \kappa_{i} \tau_{i}^{-1}\gamma_{i}^{-1}, \quad \eta_{i+1}^{-1} \lambda_{i}^{-1} \lambda_{i+1} \eta_{i} \gamma_{i}^{-1}\tau_{i}^{-1}.
\end{align*}

\item 
\begin{enumerate}[(i)]
\item If $n$ is even,
\begin{align*}
&\beta_{2}\alpha_{3}\cdots \beta_{n-2}\alpha_{n-1}\beta_{n-1}\alpha_{n-2}\cdots \beta_{3}\alpha_{2}\beta_{1} \lambda_{1}^{-1}\eta_{1}^{-1}, \quad \beta_{1}\alpha_{2}\cdots \alpha_{n-2}\beta_{n-1}\alpha_{n-1}\beta_{n-2}\cdots \alpha_{3}\beta_{2} \theta_{1}^{-1}\\
&\tau_{2}\gamma_{3}\cdots \tau_{n-2}\gamma_{n-1}\tau_{n-1}\gamma_{n-2}\cdots \tau_{3}\gamma_{2}\tau_{1}  \theta_{1}^{-1}, \quad \tau_{1}\gamma_{2}\cdots \gamma_{n-2}\tau_{n-1}\gamma_{n-1}\tau_{n-2}\cdots \gamma_{3}\tau_{2} \eta_{1}^{-1} \lambda_{1}^{-1}.
\end{align*}
\item If $n$ is odd,
\begin{align*}
&\beta_{2}\alpha_{3}\cdots \alpha_{n-2}\beta_{n-1}\alpha_{n-1}\beta_{n-2}\cdots \beta_{3}\alpha_{2}\beta_{1} \lambda_{1}^{-1}\eta_{1}^{-1}, \quad \beta_{1}\alpha_{2}\cdots \beta_{n-2}\alpha_{n-1}\beta_{n-1}\alpha_{n-2}\cdots \alpha_{3}\beta_{2} \theta_{1}^{-1}\\
&\tau_{2}\gamma_{3}\cdots \gamma_{n-2}\tau_{n-1}\gamma_{n-1}\tau_{n-2}\cdots \tau_{3}\gamma_{2}\tau_{1}  \theta_{1}^{-1}, \quad \tau_{1}\gamma_{2}\cdots \tau_{n-2}\gamma_{n-1}\tau_{n-1}\gamma_{n-2}\cdots \gamma_{3}\tau_{2} \eta_{1}^{-1} \lambda_{1}^{-1}.
\end{align*}
\end{enumerate}
\end{enumerate}
\end{enumerate}
\end{prop}

\begin{rem}
The above presentation may be used to show that $\Gamma_2(B_n(\rp))$ is perfect for $n\geq 5$.
\end{rem}

\begin{proof}
Taking the standard presentation of $B_n(\rp)$ given by \repr{present}, and the set
$\brak{1, \sigma_1, \sigma_1\rho_{1},\sigma_1\rho_{1}\sigma_{1}}$ as a Schreier transversal, we apply the Reide\-meister-Schreier rewriting process to the short exact sequence~\reqref{ses}. In this way, a generating set for $\Gamma_2(B_n(\rp))$ is that given in the statement of the proposition. We record the following equalities for later use:
\begin{equation}\label{eq:conjalgam}
\left\{
\begin{aligned}
\sigma_{1} \rho_{1} \alpha_{i} \rho_{1}^{-1} \sigma_{1}^{-1}&= \gamma_{i}\\
\sigma_{1} \rho_{1} \beta_{i} \rho_{1}^{-1} \sigma_{1}^{-1}&= \tau_{i}\\
\sigma_{1} \rho_{1} \tau_{i} \rho_{1}^{-1} \sigma_{1}^{-1}&= \sigma_{1} \rho_{1}
\sigma_{1} \rho_{1} \sigma_{1}\sigma_{i}\rho_{1}^{-1} \sigma_{1}^{-1} \rho_{1}^{-1} \sigma_{1}^{-1}=\lambda_{1}\beta_{i}\lambda_{1}^{-1}\\
\sigma_{1} \rho_{1} \gamma_{i} \rho_{1}^{-1} \sigma_{1}^{-1}&= \sigma_{1} \rho_{1}
\sigma_{1} \rho_{1} \sigma_{i}\sigma_{1}^{-1}\rho_{1}^{-1} \sigma_{1}^{-1} \rho_{1}^{-1} \sigma_{1}^{-1}=\lambda_{1}\alpha_{i}\lambda_{1}^{-1}\\
\sigma_{1} \rho_{1} \lambda_{i} \rho_{1}^{-1} \sigma_{1}^{-1}&= \sigma_{1} \rho_{1} \sigma_{1}\rho_{1} \sigma_{1}\rho_{i} \rho_{1}^{-1} \sigma_{1}^{-1}= \lambda_{1} \kappa_{i}\\
\sigma_{1} \rho_{1} \eta_{i} \rho_{1}^{-1} \sigma_{1}^{-1}&= \sigma_{1} \rho_{1} \rho_{i}\sigma_{1}^{-1}\rho_{1}^{-1}\sigma_{1}^{-1} \rho_{1}^{-1} \sigma_{1}^{-1}= \theta_{i}\lambda_{1}^{-1} \\
\sigma_{1} \rho_{1} \kappa_{i} \rho_{1}^{-1} \sigma_{1}^{-1}&= \sigma_{1} \rho_{1} 
\sigma_{1} \rho_{i}\rho_{1}^{-1}\sigma_{1}^{-1} \rho_{1}^{-1} \sigma_{1}^{-1} = \lambda_{i}\lambda_{1}^{-1}\\
\sigma_{1} \rho_{1} \theta_{i} \rho_{1}^{-1} \sigma_{1}^{-1}&= \sigma_{1} \rho_{1} 
\rho_{1} \rho_{i}\sigma_{1}^{-1}\rho_{1}^{-1}\sigma_{1}^{-1} = \lambda_{1}\eta_{i}.
\end{aligned}\right.
\end{equation}
We now determine the relations of $\Gamma_2(B_n(\rp))$ in terms of the given generating set. As we mentioned, we also set $\alpha_{1}=\kappa_{1}=\gamma_{1}=1$.  For all $1\leq i,j\leq n-1$, $\lvert i-j\rvert \geq 2$, the relator $\sigma_{i}\sigma_{j}\sigma_{i}^{-1}\sigma_{j}^{-1}$ gives rise to the following four relators, one for each element of the Schreier transversal:
\begin{align*}
1=&\sigma_{i}\sigma_{j}\sigma_{i}^{-1}\sigma_{j}^{-1}= \sigma_{i}\sigma_{1}^{-1} \ldotp \sigma_{1}\sigma_{j} \ldotp \sigma_{i}^{-1}\sigma_{1}^{-1} \ldotp \sigma_{1}\sigma_{j}^{-1}= \alpha_{i} \beta_{j} \beta_{i}^{-1}\alpha_{j}^{-1}\\
1=&\sigma_{1}\ldotp \sigma_{i}\sigma_{j}\sigma_{i}^{-1}\sigma_{j}^{-1} \ldotp \sigma_{1}^{-1}= \sigma_{1} \sigma_{i} \ldotp \sigma_{j}\sigma_{1}^{-1}\ldotp \sigma_{1} \sigma_{i}^{-1} \ldotp \sigma_{j}^{-1}\sigma_{1}^{-1}= \beta_{i} \alpha_{j} \alpha_{i}^{-1}\beta_{j}^{-1}\\
1=&\sigma_{1}\rho_{1}\ldotp \sigma_{i}\sigma_{j}\sigma_{i}^{-1}\sigma_{j}^{-1} \ldotp \rho_{1}^{-1}\sigma_{1}^{-1}= \sigma_{1} \rho_{1} \alpha_{i} \beta_{j} \beta_{i}^{-1}\alpha_{j}^{-1}\rho_{1}^{-1}\sigma_{1}^{-1}
= \gamma_{i} \tau_{j} \tau_{i}^{-1}\gamma_{j}^{-1}\\
1=&\sigma_{1}\rho_{1}\sigma_{1}\ldotp \sigma_{i}\sigma_{j}\sigma_{i}^{-1}\sigma_{j}^{-1} \ldotp \sigma_{1}^{-1}\rho_{1}^{-1}\sigma_{1}^{-1}= 
\sigma_{1} \rho_{1} \beta_{i} \alpha_{j} \alpha_{i}^{-1}\beta_{j}^{-1} \rho_{1}^{-1} \sigma_{1}^{-1} 
= \tau_{i} \gamma_{j} \gamma_{i}^{-1}\tau_{j}^{-1}.
\end{align*}
In the third and fourth equations, we have used \req{conjalgam}. Similarly, from the relator $\sigma_{i}\sigma_{i+1}\sigma_{i}\sigma_{i+1}^{-1}\sigma_{i}^{-1}\sigma_{i+1}^{-1}$, for all $1\leq i\leq n-2$, we obtain:
\begin{align*}
1=& \sigma_{i}\sigma_{i+1}\sigma_{i}\sigma_{i+1}^{-1}\sigma_{i}^{-1}\sigma_{i+1}^{-1}= \sigma_{i}\sigma_{1}^{-1}\ldotp \sigma_{1}\sigma_{i+1} \ldotp \sigma_{i}\sigma_{1}^{-1} \ldotp \sigma_{1}\sigma_{i+1}^{-1} \ldotp \sigma_{i}^{-1}\sigma_{1}^{-1} \ldotp \sigma_{1}\sigma_{i+1}^{-1}\\
=& \alpha_{i}\beta_{i+1}\alpha_{i}\alpha_{i+1}^{-1}\beta_{i}^{-1}\alpha_{i+1}^{-1}\\
1=& \sigma_{1}\ldotp \sigma_{i}\sigma_{i+1}\sigma_{i}\sigma_{i+1}^{-1}\sigma_{i}^{-1}\sigma_{i+1}^{-1} \ldotp \sigma_{1}^{-1}= 
\sigma_{1} \sigma_{i} \ldotp \sigma_{i+1}\sigma_{1}^{-1} \ldotp \sigma_{1}\sigma_{i} \ldotp \sigma_{i+1}^{-1} \sigma_{1}^{-1} \ldotp \sigma_{1}\sigma_{i}^{-1} \ldotp \sigma_{i+1}^{-1} \sigma_{1}^{-1}\\
=& \beta_{i} \alpha_{i+1} \beta_{i} \beta_{i+1}^{-1} \alpha_{i}^{-1} \beta_{i+1}^{-1}\\
1=& \sigma_{1}\rho_{1}\ldotp \sigma_{i}\sigma_{i+1}\sigma_{i}\sigma_{i+1}^{-1}\sigma_{i}^{-1}\sigma_{i+1}^{-1} \ldotp \rho_{1}^{-1}\sigma_{1}^{-1}= 
\sigma_{1} \rho_{1} \alpha_{i}\beta_{i+1}\alpha_{i}\alpha_{i+1}^{-1}\beta_{i}^{-1}\alpha_{i+1}^{-1}
\rho_{1}^{-1}\sigma_{1}^{-1}\\
=& \gamma_{i} \tau_{i+1} \gamma_{i} \gamma_{i+1}^{-1} \tau_{i}^{-1} \gamma_{i+1}^{-1}\\
1=& \sigma_{1}\rho_{1}\sigma_{1}\ldotp \sigma_{i}\sigma_{i+1}\sigma_{i}\sigma_{i+1}^{-1}\sigma_{i}^{-1}\sigma_{i+1}^{-1} \ldotp \sigma_{1}^{-1} \rho_{1}^{-1}\sigma_{1}^{-1}= 
\sigma_{1} \rho_{1} \beta_{i} \alpha_{i+1} \beta_{i} \beta_{i+1}^{-1} \alpha_{i}^{-1} \beta_{i+1}^{-1}
\rho_{1}^{-1}\sigma_{1}^{-1}\\
=& \tau_{i} \gamma_{i+1} \tau_{i} \tau_{i+1}^{-1} \gamma_{i}^{-1} \tau_{i+1}^{-1}.
\end{align*}
For all $1\leq i\leq n-1$ and $1\leq j\leq n$, $j\neq i,i+1$, the relator $\sigma_{i}\rho_{j}\sigma_{i}^{-1}\rho_{j}^{-1}$ yields:
\begin{align*}
1&= \sigma_{i}\rho_{j}\sigma_{i}^{-1}\rho_{j}^{-1}= \sigma_{i}\sigma_{1}^{-1}\ldotp \sigma_{1} \rho_{j} \rho_{1}^{-1} \sigma_{1}^{-1}\ldotp \sigma_{1} \rho_{1} \sigma_{i}^{-1} \sigma_{1}^{-1} \rho_{1}^{-1}\sigma_{1}^{-1}\ldotp \sigma_{1}  \rho_{1} \sigma_{1} \rho_{j}^{-1} = \alpha_{i}\kappa_{j}\tau_{i}^{-1}\eta_{j}^{-1}\\
1&= \sigma_{1}\ldotp \sigma_{i}\rho_{j}\sigma_{i}^{-1}\rho_{j}^{-1}\ldotp \sigma_{1}^{-1}= \sigma_{1} \sigma_{i} \ldotp \rho_{j} \sigma_{1}^{-1} \rho_{1}^{-1}\sigma_{1}^{-1}   \ldotp \sigma_{1} \rho_{1} \sigma_{1}\sigma_{i}^{-1} \rho_{1}^{-1}\sigma_{1}^{-1}\ldotp
\sigma_{1}  \rho_{1} \rho_{j}^{-1} \sigma_{1}^{-1}\\
& = \beta_{i}\eta_{j}\gamma_{i}^{-1}\kappa_{j}^{-1}\\
1&= \sigma_{1}\rho_{1}\ldotp \sigma_{i}\rho_{j}\sigma_{i}^{-1}\rho_{j}^{-1}\ldotp \rho_{1}^{-1} \sigma_{1}^{-1}= 
\sigma_{1} \rho_{1} \alpha_{i}\kappa_{j}\tau_{i}^{-1}\kappa_{j}^{-1} \rho_{1}^{-1}\sigma_{1}^{-1}= \gamma_{i}\lambda_{j}\beta_{i}^{-1}\theta_{j}^{-1}\\
1&= \sigma_{1}\rho_{1}\sigma_{1} \ldotp \sigma_{i}\rho_{j}\sigma_{i}^{-1}\rho_{j}^{-1}\ldotp \sigma_{1}^{-1} \rho_{1}^{-1} \sigma_{1}^{-1}= 
\sigma_{1} \rho_{1} \beta_{i}\eta_{j}\gamma_{i}^{-1}\kappa_{j}^{-1}\rho_{1}^{-1}\sigma_{1}^{-1}= \tau_{i}\theta_{j}\alpha_{i}^{-1}\lambda_{j}^{-1}.
\end{align*}
For all $1\leq i\leq n-1$, the relator $\sigma_{i}^{-1} \rho_{i}\sigma_{i}^{-1} \rho_{i+1}^{-1}$ gives rise to:
\begin{align*}
1&=\sigma_{i}^{-1} \rho_{i}\sigma_{i}^{-1} \rho_{i+1}^{-1}= \sigma_{i}^{-1}\sigma_{1}^{-1} \ldotp \sigma_{1} \rho_{i}\rho_{1}^{-1} \sigma_{1}^{-1} \ldotp \sigma_{1} \rho_{1} \sigma_{i}^{-1}\sigma_{1}^{-1} \rho_{1}^{-1}\sigma_{1}^{-1}\ldotp \sigma_{1}  \rho_{1} \sigma_{1} \rho_{i+1}^{-1}\\
&= \beta_{i}^{-1} \kappa_{i} \tau_{i}^{-1} \eta_{i+1}^{-1}\\
1&=\sigma_{1}\ldotp \sigma_{i}^{-1} \rho_{i}\sigma_{i}^{-1} \rho_{i+1}^{-1}\ldotp \sigma_{1}^{-1}= 
\sigma_{1}\sigma_{i}^{-1}\ldotp \rho_{i}\sigma_{1}^{-1} \rho_{1}^{-1} \sigma_{1}^{-1} \ldotp 
\sigma_{1} \rho_{1} \sigma_{1} \sigma_{i}^{-1} \rho_{1}^{-1}\sigma_{1}^{-1}\ldotp \sigma_{1}  \rho_{1} \rho_{i+1}^{-1}\sigma_{1}^{-1}\\
& = \alpha_{i}^{-1} \eta_{i} \gamma_{i}^{-1} \kappa_{i+1}^{-1} \\
1&=\sigma_{1}\rho_{1}\ldotp \sigma_{i}^{-1} \rho_{i}\sigma_{i}^{-1} \rho_{i+1}^{-1}\ldotp \rho_{1}^{-1} \sigma_{1}^{-1}= 
\sigma_{1} \rho_{1}  \beta_{i}^{-1} \kappa_{i} \tau_{i}^{-1} \eta_{i+1}^{-1} \rho_{1}^{-1} \sigma_{i}^{-1}=\tau_{i}^{-1} \lambda_{i} \beta_{i}^{-1} \theta_{i+1}^{-1} \\
1&=\sigma_{1}\rho_{1}\sigma_{1} \ldotp \sigma_{i}^{-1} \rho_{i}\sigma_{i}^{-1} \rho_{i+1}^{-1}\ldotp \sigma_{1}^{-1} \rho_{1}^{-1} \sigma_{1}^{-1}= 
\sigma_{1} \rho_{1} \alpha_{i}^{-1} \eta_{i} \gamma_{i}^{-1} \kappa_{i+1}^{-1}  \rho_{1}^{-1} \sigma_{i}^{-1}= \gamma_{i}^{-1} \theta_{i} \alpha_{i}^{-1} \lambda_{i+1}^{-1}. 
\end{align*}
From the relator $\rho_{i+1}^{-1}\rho_{i}^{-1}\rho_{i+1}\rho_{i}\sigma_{i}^{-2}$, for all $1\leq i\leq n-1$, we obtain:
\begin{align*}
1=&\rho_{i+1}^{-1}\rho_{i}^{-1}\rho_{i+1}\rho_{i}\sigma_{i}^{-2}= \rho_{i+1}^{-1}\sigma_{1}^{-1} \rho_{1}^{-1}\sigma_{1}^{-1}\ldotp \sigma_{1}\rho_{1}\sigma_{1}\rho_{i}^{-1} \ldotp 
\rho_{i+1} \rho_{1}^{-1} \rho_{1}^{-1}\sigma_{1}^{-1}\ldotp \sigma_{1}\rho_{1}\sigma_{1}\rho_{i}\ldotp\\
&\sigma_{i}^{-1}\sigma_{1}^{-1} \ldotp \sigma_{1}\sigma_{i}^{-1} =\lambda_{i+1}^{-1} \eta_{i}^{-1} \eta_{i+1} \lambda_{i}\beta_{i}^{-1} \alpha_{i}^{-1}\\
1=&\sigma_{1} \ldotp \rho_{i+1}^{-1}\rho_{i}^{-1}\rho_{i+1}\rho_{i}\sigma_{i}^{-2} \ldotp \sigma_{1}^{-1}= 
\sigma_{1} \rho_{i+1}^{-1}\rho_{1}^{-1} \sigma_{1}^{-1}\ldotp \sigma_{1}\rho_{1} \rho_{i}^{-1} \sigma_{1}^{-1} \ldotp 
\sigma_{1} \rho_{i+1} \rho_{1}^{-1} \sigma_{1}^{-1}\ldotp\\
& \sigma_{1}\rho_{1} \rho_{i}\sigma_{1}^{-1}\ldotp
\sigma_{1} \sigma_{i}^{-1} \ldotp \sigma_{i}^{-1}\sigma_{1}^{-1}=\theta_{i+1}^{-1} \kappa_{i}^{-1} \kappa_{i+1} \theta_{i} \alpha_{i}^{-1}\beta_{i}^{-1}\\
1=&\sigma_{1} \rho_{1}\ldotp \rho_{i+1}^{-1}\rho_{i}^{-1}\rho_{i+1}\rho_{i}\sigma_{i}^{-2} \ldotp \rho_{1}^{-1}\sigma_{1}^{-1}= 
\sigma_{1} \rho_{1} \lambda_{i+1}^{-1} \eta_{i}^{-1} \eta_{i+1} \lambda_{i}\beta_{i}^{-1} \alpha_{i}^{-1} \rho_{1}^{-1} \sigma_{1}^{-1}\\
=& \kappa_{i+1}^{-1} \theta_{i}^{-1} \theta_{i+1} \kappa_{i} \tau_{i}^{-1}\gamma_{i}^{-1}\\
1=&\sigma_{1} \rho_{1}\sigma_{1}\ldotp \rho_{i+1}^{-1}\rho_{i}^{-1}\rho_{i+1}\rho_{i}\sigma_{i}^{-2} \ldotp \sigma_{1}^{-1}\rho_{1}^{-1}\sigma_{1}^{-1}= 
\sigma_{1} \rho_{1} \theta_{i+1}^{-1} \kappa_{i}^{-1} \kappa_{i+1} \theta_{i} \alpha_{i}^{-1}\beta_{i}^{-1} \rho_{1}^{-1} \sigma_{1}^{-1}\\
=&\eta_{i+1}^{-1} \lambda_{i}^{-1} \lambda_{i+1} \eta_{i} \gamma_{i}^{-1}\tau_{i}^{-1}.
\end{align*}
Finally we come to the surface relator $\sigma_{1}\cdots \sigma_{n-1}\sigma_{n-1}\cdots \sigma_{1}\rho_{1}^{-2}$. We deal with the cases $n$ even and odd separately.
\begin{enumerate}[(a)]
\item \underline{$n$ even.} We have:
\begin{align*}
1=&\sigma_{1}\cdots \sigma_{n-1}\sigma_{n-1}\cdots \sigma_{1}\rho_{1}^{-2}= \sigma_{1}\sigma_{2} \ldotp \sigma_{3}\sigma_{1}^{-1} \cdots \sigma_{1}\sigma_{n-2} \ldotp \sigma_{n-1}\sigma_{1}^{-1} \ldotp \sigma_{1} \sigma_{n-1} \ldotp \sigma_{n-2}\sigma_{1}^{-1} \cdots\\
& \sigma_{1} \sigma_{3} \ldotp \sigma_{2}\sigma_{1}^{-1}\ldotp \sigma_{1}^2\ldotp \rho_{1}^{-1} \sigma_{1}^{-1} \rho_{1}^{-1}\sigma_{1}^{-1} \ldotp \sigma_{1} \rho_{1} \sigma_{1} \rho_{1}^{-1}
\\ 
=& \beta_{2}\alpha_{3}\cdots \beta_{n-2}\alpha_{n-1}\beta_{n-1}\alpha_{n-2}\cdots \beta_{3}\alpha_{2}\beta_{1} \lambda_{1}^{-1}\eta_{1}^{-1}\\
1 =& \sigma_{1}\ldotp \sigma_{1}\cdots \sigma_{n-1}\sigma_{n-1}\cdots \sigma_{1}\rho_{1}^{-2} \ldotp \sigma_{1}^{-1} = \sigma_{1}^2 \ldotp \sigma_{2} \sigma_{1}^{-1} \ldotp \sigma_{1}\sigma_{3} \cdots \sigma_{n-2} \sigma_{1}^{-1}\ldotp \sigma_{1}\sigma_{n-1} \ldotp \sigma_{n-1} \sigma_{1}^{-1}  \ldotp\\ & \sigma_{1}\sigma_{n-2} \cdots
\sigma_{3} \sigma_{1}^{-1} \ldotp \sigma_{1}\sigma_{2}\ldotp \sigma_{1} \rho_{1}^{-2} \sigma_{1}^{-1} 
= \beta_{1}\alpha_{2}\cdots \alpha_{n-2}\beta_{n-1}\alpha_{n-1}\beta_{n-2}\cdots \alpha_{3}\beta_{2} \theta_{1}^{-1}\\
1 =& \sigma_{1} \rho_{1} \ldotp \sigma_{1}\cdots \sigma_{n-1}\sigma_{n-1}\cdots \sigma_{1}\rho_{1}^{-2} \ldotp \rho_{1}^{-1} \sigma_{1}^{-1}\\
=& \sigma_{1} \rho_{1} \ldotp \beta_{2}\alpha_{3}\cdots \beta_{n-2}\alpha_{n-1}\beta_{n-1}\alpha_{n-2}\cdots \beta_{3}\alpha_{2}\beta_{1} \lambda_{1}^{-1}\eta_{1}^{-1} \ldotp \rho_{1}^{-1} \sigma_{1}^{-1}\\
=&\tau_{2}\gamma_{3}\cdots \tau_{n-2}\gamma_{n-1}\tau_{n-1}\gamma_{n-2}\cdots \tau_{3}\gamma_{2}\tau_{1}  \theta_{1}^{-1}\\
1 =& \sigma_{1} \rho_{1}\sigma_{1}\ldotp \sigma_{1}\cdots \sigma_{n-1}\sigma_{n-1}\cdots \sigma_{1}\rho_{1}^{-2} \ldotp \sigma_{1}^{-1} \rho_{1}^{-1} \sigma_{1}^{-1}\\
=& \sigma_{1} \rho_{1}\ldotp
\beta_{1}\alpha_{2}\cdots \alpha_{n-2}\beta_{n-1}\alpha_{n-1}\beta_{n-2}\cdots \alpha_{3}\beta_{2} \theta_{1}^{-1} \ldotp
\rho_{1}^{-1} \sigma_{1}^{-1}
\\
=& \tau_{1}\gamma_{2}\cdots \gamma_{n-2}\tau_{n-1}\gamma_{n-1}\tau_{n-2}\cdots \gamma_{3}\tau_{2} \eta_{1}^{-1} \lambda_{1}^{-1}.
\end{align*}

\item \underline{$n$ odd.} We have:
\begin{align*}
1=&\sigma_{1}\cdots \sigma_{n-1}\sigma_{n-1}\cdots \sigma_{1}\rho_{1}^{-2}= \sigma_{1}\sigma_{2} \ldotp \sigma_{3}\sigma_{1}^{-1} \cdots \sigma_{n-2}\sigma_{1}^{-1} \ldotp \sigma_{1}\sigma_{n-1} \ldotp  \sigma_{n-1} \sigma_{1}^{-1} \ldotp \sigma_{1}\sigma_{n-2} \cdots\\
& \sigma_{1} \sigma_{3} \ldotp \sigma_{2}\sigma_{1}^{-1}\ldotp \sigma_{1}^2\ldotp \rho_{1}^{-1} \sigma_{1}^{-1} \rho_{1}^{-1}\sigma_{1}^{-1} \ldotp \sigma_{1} \rho_{1} \sigma_{1} \rho_{1}^{-1}\\
=& \beta_{2}\alpha_{3}\cdots \alpha_{n-2}\beta_{n-1}\alpha_{n-1}\beta_{n-2}\cdots \beta_{3}\alpha_{2}\beta_{1} \lambda_{1}^{-1}\eta_{1}^{-1}\\
1 =& \sigma_{1}\ldotp \sigma_{1}\cdots \sigma_{n-1}\sigma_{n-1}\cdots \sigma_{1}\rho_{1}^{-2} \ldotp \sigma_{1}^{-1} = \sigma_{1}^2 \ldotp \sigma_{2} \sigma_{1}^{-1} \ldotp \sigma_{1}\sigma_{3} \cdots \sigma_{1} \sigma_{n-2} \ldotp \sigma_{n-1} \sigma_{1}^{-1} \ldotp \sigma_{1} \sigma_{n-1}  \ldotp\\
& \sigma_{n-2}\sigma_{1}^{-1} \cdots
\sigma_{3} \sigma_{1}^{-1} \ldotp \sigma_{1}\sigma_{2}\ldotp \sigma_{1} \rho_{1}^{-2} \sigma_{1}^{-1} 
= \beta_{1}\alpha_{2}\cdots \beta_{n-2}\alpha_{n-1}\beta_{n-1}\alpha_{n-2}\cdots \alpha_{3}\beta_{2} \theta_{1}^{-1}\\
1 =& \sigma_{1} \rho_{1} \ldotp \sigma_{1}\cdots \sigma_{n-1}\sigma_{n-1}\cdots \sigma_{1}\rho_{1}^{-2} \ldotp \rho_{1}^{-1} \sigma_{1}^{-1}\\
=& \sigma_{1} \rho_{1} \ldotp \beta_{2}\alpha_{3}\cdots \alpha_{n-2}\beta_{n-1}\alpha_{n-1}\beta_{n-2}\cdots \beta_{3}\alpha_{2}\beta_{1} \lambda_{1}^{-1}\eta_{1}^{-1} \ldotp \rho_{1}^{-1} \sigma_{1}^{-1}\\
=& \tau_{2}\gamma_{3}\cdots \gamma_{n-2}\tau_{n-1}\gamma_{n-1}\tau_{n-2}\cdots \tau_{3}\gamma_{2}\tau_{1} \lambda_{1}^{-1}\lambda_{1} \theta_{1}^{-1}\\
1 =& \sigma_{1} \rho_{1}\sigma_{1}\ldotp \sigma_{1}\cdots \sigma_{n-1}\sigma_{n-1}\cdots \sigma_{1}\rho_{1}^{-2} \ldotp \sigma_{1}^{-1} \rho_{1}^{-1} \sigma_{1}^{-1}\\
=& \sigma_{1} \rho_{1}\ldotp
\beta_{1}\alpha_{2}\cdots \beta_{n-2}\alpha_{n-1}\beta_{n-1}\alpha_{n-2}\cdots \alpha_{3}\beta_{2} \theta_{1}^{-1} \ldotp
\rho_{1}^{-1} \sigma_{1}^{-1}
\\
=& \tau_{1}\gamma_{2}\cdots \tau_{n-2}\gamma_{n-1}\tau_{n-1}\gamma_{n-2}\cdots \gamma_{3}\tau_{2} \eta_{1}^{-1} \lambda_{1}^{-1}.
\end{align*}
This completes the proof of the proposition.\qedhere
\end{enumerate}
\end{proof}

\end{document}